\documentclass[a4paper,11pt,oneside]{article}

\usepackage{amsmath,amsthm,amsfonts,amssymb,amscd}
\usepackage{graphicx}
\usepackage{mathtools}
\usepackage{hyperref}
\usepackage{mathrsfs}
\usepackage{float}
\usepackage{csvsimple}
\usepackage{enumerate}
\usepackage{booktabs}
\usepackage{pgfplots}
\usepackage[noend]{algpseudocode}
\usepackage{graphicx,bm,xcolor}
\usepackage{algorithm}
\usepackage{pdfpages}
\usepackage{algpseudocode}
\usepackage{stmaryrd}
\usepackage{hyperref}
\usetikzlibrary{arrows}
\usepackage{caption}
\usepackage[T1]{fontenc}
\usepackage{t1enc}
\usepackage[ngerman,english]{babel}
\usepackage{verbatim} 

\usepackage{url}

\usepackage{units}  
\usepackage{siunitx}

\newcommand{\triple}{\ensuremath{| \! | \! |}} 
\newcommand{\proofbegin}{\noindent\textsc{Proof }}
\newcommand{\proofend}{\hfill$\square$}

\usepackage{tabularx}
\usepackage[font={footnotesize},%
  labelsep=colon]{caption}
\usepackage{tabu}
\usepackage{tikz}
\tikzset{every path/.append style={line width=1pt}}

\usepackage{authblk}

\usepackage[top=2cm, bottom=2cm, left=2cm, right=2cm]{geometry}

\theoremstyle{plain} 
\newtheorem{theorem}{Theorem}[section]
\newtheorem{Proposition}{Proposition}[section]
\newtheorem{lemma}[theorem]{Lemma}

\newtheorem{Remark}[theorem]{Remark}

\newtheorem{form}[theorem]{Formulation}

\theoremstyle{definition} %

\theoremstyle{remark} %

\begin{document}

\title{A phase-field model for fractures in incompressible solids}
\author[1]{Katrin Mang}
\author[1]{Thomas Wick}
\author[2]{Winnifried Wollner}

\affil[1]{Leibniz Universit\"at Hannover, Institut f\"ur Angewandte
  Mathematik, AG Wissenschaftliches Rechnen, Welfengarten 1, 30167 Hannover, Germany}
\affil[2]{Department of Mathematics, Technische Universit{\"a}t Darmstadt,
Dolivostrasse 15, 64293 Darmstadt, Germany}

\date{}

\maketitle
	
\begin{abstract}
Within this work, we develop a phase-field description for simulating 
fractures in incompressible materials.
Standard formulations are subject to volume-locking when 
the solid is (nearly) incompressible. 
We propose an approach that builds on a mixed form of the displacement
equation with two unknowns: 
a displacement field and a hydro-static pressure variable. Corresponding
function spaces have to be chosen properly. On the discrete level, stable
Taylor-Hood elements are employed for the displacement-pressure system. 
Two additional variables describe the phase-field solution and the 
crack irreversibility constraint. Therefore, the final 
system contains four variables: displacements, pressure, phase-field, and 
a Lagrange multiplier. The resulting discrete system is nonlinear and 
solved monolithically with a Newton-type method.
Our proposed model is demonstrated by means of several numerical studies based on two numerical tests.
First, different finite element choices are compared in order to investigate the influence 
of higher-order elements in the proposed settings.
Further, numerical results including spatial mesh refinement studies and variations
in Poisson's ratio approaching the incompressible limit, are presented.
\end{abstract}

\section{Introduction}\label{sec_intro}
Currently, crack propagation is one of the major research topics in mechanical, energy, and environmental engineering. A well-established variational approach for Griffith's \cite{griffith1920phenomena} 
quasi-static brittle fracture was introduced by Francfort
and Marigo \cite{FraMar98}. 
Since then, the method was applied in numerous different studies in calculus
of variations, numerical analysis, and engineering.
Miehe et al. \cite{MieWelHof10a} introduced the name phase-field modeling for this
variational approach.

To the best of our knowledge, in all published studies, it is assumed that the material has a Poisson
ratio $\nu$ much smaller than $0.5$. It implies dealing with compressible solids. Thus, the objective of this work is handling fractures in (nearly) incompressible materials.
Incompressible solids are important for various sciences, i.e. engineering \cite{holzapfel2002nonlinear,taylor2011isogeometric} or medicine \cite{holzapfel1996large,schroder2005variational}.
One industrial example is the design of reliable rubber products
\cite{kubo2017velocity}, 
which gives a clear argument to improve the understanding of the 
mechanical properties and fracture mechanism of incompressible materials.
The special properties of incompressible solids make it
challenging to successfully explain or reproduce
crack propagation in rubbers via numerical simulations \cite{vsuvstarivc2014sensitivity}.
The novel aspect of this work is to investigate 
cases with Poisson's ratio $\nu$ approximating $0.5$, which relates to incompressible
materials such as for instance rubber. 

The ratio of the Lam\'{e} coefficient $\lambda$ to the Poisson ratio $\nu$ and the Lam\'{e} coefficient
$\mu$ is given by
\begin{align*}
\lambda=\frac{2 \nu \mu}{1-2 \mu}.
\end{align*}
If $\nu$ tends to $0.5$, the parameter $\lambda$ increases and becomes much
larger than $\mu$. This situation is well-known in solid mechanical simulations
as so-called Poisson or volume-locking \cite{babuvska1992locking}.

One possibility to avoid these effects is a Discontinuous Galerkin (DG) method, e.g., \cite{cockburn2000development}. 
Whiler \cite{wihler2006locking} used the DG method for linear elasticity problems, 
Hansbo et al. \cite{hansbo2002discontinuous} studied in particular incompressible and nearly incompressible elasticity problems.
Another possibility is to split the displacement equation into a mixed system, see for
instance Braess \cite{braess2007finite}. The major problem of locking is that 
C\'ea's lemma delivers errors which can be significantly larger than the approximation error. 

In this work, we concentrate on a mixed problem formulation. 
The displacement equation is split into a modified displacement equation
for computing $u$ and an equation determining a hydro-static pressure $p$.
To ensure
stability, an inf-sup condition \cite{girault2012finite} must be fulfilled. This 
means that the respective sets in the mixed formulation must be carefully
chosen. Then, this condition carries over to the discrete mixed system.
Here, the discrete space for $u$ must be larger than the space for $p$. 
As finite element approach, we make use of the Taylor-Hood element with
biquadratic shape functions ($Q_2$) for the displacement field and bilinear shape
functions ($Q_1$) for the pressure. Indeed, $Q_2 Q_1$ elements fulfill the discrete
inf-sup condition.

To account for crack irreversibility (the crack cannot heal), 
the phase-field fracture formulation turns to a variational inequality.
To treat the inequality constraint we employ a Lagrange multiplier, see e.g.,
\cite{ito2008lagrange,rockafellar1993lagrange}.

The main contributions of this work are:
\begin{itemize}
 \item Formulating a quasi-static phase-field fracture model for incompressible solids;
 \item Discretizing the new model with stable finite elements;
 \item Substantiating the new model with appropriate numerical tests including studies with varying the spatial discretization parameter $h$ and Poisson's ratio $\nu$.
\end{itemize}

The outline of this paper is as follows: In Section \ref{sec_eq}, the notation and the basic system of equations are presented.
Afterward, the new model formulation is proposed in Section \ref{NewModel} and numerical tests, based on the new model, are presented in the subsequent Section \ref{Results}.
Two well-known mechanical tests are evaluated, in particular, to investigate the consistency of the phase-field fracture model in a mixed form. 
To strengthen our findings, results using the new model formulation are compared with the common quasi-static phase-field model and the standardly used finite element approach.
Furthermore, numerical solutions on finer meshes and with different Poisson ratios up to a nearly incompressible setting are exploited. 
Section \ref{Conclusion} summarizes the content of this work.


\section{Notation}\label{sec_eq}
This section covers the basic notation, the function spaces, the required variables and the standard problem formulation for a phase-field approach.\\
We emanate from a two-dimensional, open and smooth domain $\Omega\subset \mathbb{R}^2$. By means of an elliptic functional developed 
by Ambrosio-Tortorelli \cite{ambrosio1992approximation,ambrosio1990approximation}, 
it exists a lower-dimensional crack $C \subset
\mathbb{R}^1 \in \Omega$. On the boundary $\partial \Omega$ we assume to have homogeneous
Dirichlet boundary conditions. Let $I$ be a loading (time) interval $(0,T)$,
where $T>0$ is the end time value.
A displacement function $u:(\Omega \times I) \to \mathbb{R}^2$ is defined on the domain $\Omega$.
Introducing the phase-field approach, the crack is approximated via a phase-field variable $\varphi: (\Omega \times I) \to [0,1]$ with $\varphi=0$ in the crack and $\varphi=1$ in the unbroken material. 
A parameter $\epsilon>0$ determines the width of a transition zone between the unbroken material and the broken material inside the crack $C$.
To handle the irreversibility constraint, we use a Lagrange multiplier 
$\tau:(\Omega\times I) \to \mathbb{R}$ with $\tau \geq 0$.\\
The Frobenius scalar product of two matrices of the same dimension is defined as $(A:B):= \sum_i \sum_j a_{ij} b_{ij}$. 
By $(a,b) := \int_\Omega a \cdot b\ dx$ for vectors $a, b$ the $L^2$ scalar-product is denoted. 
For tensor-valued functions $A$ and $B$ of the same dimension it holds $(A,B) := \int_\Omega A : B\ dx$. The expression $\|.\|_n$ denominates the Sobolev-norm of order $n$. 
The maximum of two values is denoted by $\max\{.,.\}$.\\

For a complete formulation of the phase-field model, further definitions are needed. A degradation function $g(\varphi)$ is defined as
\begin{align*}
g(\varphi):=(1-\kappa)\varphi^2 + \kappa,
\end{align*}
with a small regularization parameter $\kappa > 0$.
The stress tensor $\sigma(u)$ is given by
\begin{align*}
\sigma(u) := 2 E_{\text{lin}}(u) + \lambda \text{tr} (E_{\text{lin}}(u)) \textbf{I}.
\end{align*} 
Next, $E_{\text{lin}}(u)$ is the linearized strain tensor:
\begin{align*}
E_{\text{lin}}(u):=\frac{1}{2} (\nabla u + \nabla u^T). 
\end{align*}
By $\textbf{I}$, the two-dimensional identity matrix is denoted.
The continuous formulation referred to Miehe et al. \cite{MieWelHof10a} is given in the following. \\
Find $u:(\Omega\times I) \to \mathbb{R}^2$ and $\varphi:(\Omega \times I) \to \mathbb{R}$ such that
\begin{align*}
-\nabla \cdot (g(\varphi) \sigma(u)) = 0\quad \text{in}\ (\Omega \times I),
\end{align*}
\begin{align}
\begin{aligned}
(1- \kappa)\varphi E_{\text{lin}}(u): E_{\text{lin}}(u)- \frac{G_c}{\epsilon}(1-\varphi)
 + \epsilon G_c \Delta \varphi \geq 0\quad \text{in}\ (\Omega \times I).\label{phi}
\end{aligned}
\end{align}
Herein, $G_c$  is the critical energy release rate.
The crack irreversibility condition is determined by
\begin{align}
 \partial_t \varphi \leq 0\quad \text{in}\ (\Omega \times I),\label{irre}
\end{align}
which has to be cautiously treated in the numerical solving. In the frame of this work, the crack irreversibility is discretized
via $\varphi^n\leq\varphi^{n-1}$ for loading increments $n$.

A Dirichlet boundary condition for the displacement function $u$ can be stated as 
\begin{align*}
  u= u_D\quad \text{on}\ (\partial\Omega\times I).
\end{align*}
To link the phase-field equation in (\ref{phi}) and the crack irreversibility constraint in (\ref{irre}), a compatibility condition is required:
\begin{align}
\begin{aligned}
\Bigl((1- \kappa)\varphi E_{\text{lin}}(u):E_{\text{lin}}(u)- \frac{G_c}{\epsilon}(1- \varphi)
+ G_c \epsilon \Delta \varphi \Bigr)\cdot (\partial_t \varphi) = 0\quad \text{in}\ (\Omega \times I).\label{comp}
\end{aligned}
\end{align}
With
\begin{align*}
\varphi(x,0)= \varphi_0\quad \text{in}\ (\Omega \times \{0\}),
\end{align*}
an initial condition is imposed to complete the problem formulation.
Based on the continuous phase-field problem, in the next section the corresponding variational problem formulation is provided. Afterwards, we present a new mixed problem formulation.


\section{A Phase-field Model for Incompressible Solids}\label{NewModel}
Within this section, a stable phase-field formulation for incompressible materials is presented.  
First, the variational formulation of the quasi-static phase-field approach is given with suitable ansatz spaces. 
\begin{Remark}
From now on, we assume to deal with a time-discretized (in mechanics:
incremental) problem in which the loading/time interval $I = (0,T)$ is
discretized
using the time points:
\[
0=t_0<t_1<t_2 < \ldots < t_n < \ldots <t_N=T.
\]
To this end, the irreversibility constraint \eqref{irre} is approximated 
as 
\[
\varphi^n\leq\varphi^{n-1}
\]
with $\varphi^n := \varphi(t_n)$ and $\varphi^{n-1} := \varphi(t_{n-1})$.
\end{Remark}


\subsection{A variational phase-field fracture formulation}\label{incompPhasefield}
We first recall an often employed variational formulation 
for quasi-static brittle fracture. To this end, 
the discretized irreversibility condition
 $\varphi^n \leq \varphi^{n-1}$
is embedded in the feasible set for the phase-field variable.
We define function spaces $\mathcal{V}:= H_0^1(\Omega)^2$,
$\mathcal{W}:=  H^1(\Omega)$, a convex subset $\mathcal{K}:=\mathcal{K}(\varphi^{n-1})\subset \mathcal{W}$ and for later purposes $\mathcal{U}:=L_2(\Omega)$ and $\mathcal{X}:=L_2(\Omega)$.


\begin{form}[Variational problem]\label{formVariational}\hfill

The weak form of the phase-field fracture problem reads as follows:\\
Find $u \in \{u_D + \mathcal{V}\}$ and $\varphi \in \mathcal{W}$ such that 
 \begin{align}
 \begin{aligned}
  2\mu(g(\varphi) E_{\text{lin}}(u),E_{\text{lin}}(w))
  + \lambda(g(\varphi)\nabla \cdot u, \nabla\cdot w)= 0\quad \forall w \in \mathcal{V},\label{uequation}
  \end{aligned}
 \end{align}
 \begin{align*}
 \begin{aligned}
  2(1- \kappa) \Bigl(\varphi E_{\text{lin}}(u):&\ E_{\text{lin}}(u),\psi-\varphi \Bigr) 
+ (1-\kappa)(\varphi \lambda \nabla \cdot u \textbf{I}:E_{\text{lin}}(u),\psi-\varphi)\\
+&\ G_c (-\frac{1}{\epsilon} (1-\varphi),\psi-\varphi)
+ G_c \epsilon(\nabla \varphi,\nabla(\psi-\varphi))\geq 0\quad \forall \psi \in \mathcal{K}. 
 \end{aligned}
 \end{align*}
\end{form}


\subsection{Investigation of incompressible solids for a decoupled system}
\label{sec_decoupled}
In this section, we highlight the problem of incompressible solids for a
simplified decoupled problem. We restrict ourselves
to the displacement equation \eqref{uequation} assuming that the phase-field 
variable is a given coefficient with sufficient regularity.

As mentioned in the introduction, it holds $\lambda \gg \mu$ for (nearly) incompressible solids.
This results in a large increase of the energy within small density changes. To simplify discussing the stability of the $u$-equation containing a large parameter $\lambda$, 
we assume homogeneous Dirichlet boundary conditions on $\partial\Omega$.

Due to the decoupling, the $u$-equation is linear which allows us to apply the 
usual properties to obtain a well-posed problem via the Lax-Milgram lemma.
We define the bilinear form 
\begin{align*}
 a_{\varphi}(u,w):= (g(\varphi) E_{\text{lin}}(u),E_{\text{lin}}(w)).
\end{align*}
For the estimates of continuity and coercivity, we emphasize that 
a coefficient $g(\varphi)$ enters. First, we assume 
$\kappa$ to be small, but constant throughout this paper.
For $\varphi = 0$ (in the fracture zone), we have
\[
g(\varphi) = \kappa.
\]
In the unbroken material, i.e. $\varphi = 1$, it holds
\[
g(\varphi) = 1.
\]
In particular, we define:
\[
\alpha := \inf_{x \in \Omega} \alpha_0 g(\varphi(x)), \quad C:= \sup_{x\in \Omega} C_0 g(\varphi(x))
\]
for given $\alpha_0 > 0$ and $C_0 > 0$. The constants $\alpha_0$ and $C_0$
arise in 
considering the coercivity and continuity of $-\nabla \cdot (g(\varphi) \sigma(u))$. We recall,
that via
\begin{align*}
\begin{aligned}
|a_{\varphi}(u,w)|\leq&\ C \|u\|_{\mathcal{V}} \|w\|_{\mathcal{V}}\quad \text{for}\ C >\ 0\quad \forall u,w \in \mathcal{V},\\[3pt]
 a_{\varphi}(w,w)\geq&\ \alpha \|w\|_{\mathcal{V}}^2 \quad\quad\quad\ \text{for}\ \alpha >\ 0\ \quad\ \forall w \in \mathcal{V},
\end{aligned}
\end{align*}
the bilinear form $a(\cdot,\cdot)$
is continuous and $H^1$-elliptic.
Notice, that for small $\kappa$ (hidden in $\alpha$ through $g(\varphi)$) the coercivity 
estimate may become critical. 
Furthermore, it holds $\alpha \leq \mu$ and $C\geq \lambda + \mu$, see, e.g., \cite{braess2007finite}. 
Consequently, we obtain 
\begin{align*}
\frac{C}{\alpha} \to \infty \quad\text{for } \lambda\to \infty.
\end{align*}
This estimate becomes even worse when $\kappa \approx 0$.
Because $\frac{C}{\alpha}$ enters into the C\'ea lemma, 
we obtain larger errors than the expected approximation errors. 
This phenomenon is called volume-locking \cite{babuvska1992locking}.

As described in the introductory part, one approach to avoid locking is a mixed problem formulation with penalty term. Using this technique, we define
\begin{align*}
 p:= \lambda \nabla \cdot u\quad\text{with}\  p \in \mathcal{U}.
\end{align*}
As we will see later, the variable $p$ describes a hydro-static pressure.
It allows to reformulate the $u$-equation into a mixed system:

Find $u \in \mathcal{V}$ and $p \in \mathcal{U}$ such that
 \begin{align}
 \begin{aligned}
2 \mu (g(\varphi) E_{\text{lin}}(u),E_{\text{lin}}(w)) + (g(\varphi) p, \nabla \cdot w) =&\ 0\quad \forall w \in \mathcal{V},\label{mixed1}\\
(g(\varphi)\nabla \cdot u,q) - \frac{1}{\lambda} (g(\varphi)p,q) =&\ 0\quad \forall q \in \mathcal{U}.
\end{aligned}
\end{align}
\begin{Remark}
The coefficient $\lambda$ arises in the denominator in the mixed formulation. For
this reason, a large $\lambda$ is less harmful.
\end{Remark}
\begin{Remark}
The previous system is a saddle-point problem with penalty term. The penalty 
parameter is nothing else than the Lam\'{e} coefficient $\lambda$.
\end{Remark}

To simplify the notation (and to adapt to the literature), 
we introduce the following bilinear forms:
\begin{align*}
 \begin{aligned}
  a_{\varphi}(u,w)=&\ (g(\varphi) E_{\text{lin}}(u),E_{\text{lin}}(w)),\\
  b_{\varphi}(w,p):=&\ (g(\varphi) \nabla \cdot w,p),\\
  c_{\varphi}(p,q):=&\ (g(\varphi)p,q).
 \end{aligned}
\end{align*}
For $c_{\varphi}(p,q)$ we define a semi-norm $|q|_c:=c_{\varphi}(q,q)^{\frac{1}{2}}$. Because $a_{\varphi}(u,w)=a_{\varphi}(w,u)$ and $c_{\varphi}(p,q)=c_{\varphi}(q,p)$, the bilinear form $a$ and $c$ are symmetric.
Then the previous system can be stated as:\\
Find $(u,p) \in (\mathcal{V} \times \mathcal{U})$ such that
 \begin{align}
   2 \mu a_{\varphi}(u,w) + b_{\varphi}(w,p) = 0 &\ \forall w \in \mathcal{V},\label{41}\\
   b_{\varphi}(u,q) - \frac{1}{\lambda} c_{\varphi}(p,q) = 0 &\ \forall q \in \mathcal{U}\label{42}.
 \end{align}
For the following, we define a
compact bilinear form summing up the single terms from before:
\begin{align*}
\begin{aligned}
 A_{\varphi}(u,p;w,q) := 2 \mu a_{\varphi}(u,w) + b_{\varphi}(w,p) + b_{\varphi}(u,q) - \frac{1}{\lambda} c_{\varphi}(p,q).
 \end{aligned}
\end{align*}
Furthermore, the natural norm for a saddle point problem is defined as
\begin{align*}
 \triple(w,q)\triple:= \|w\|_{\mathcal{V}} + \|q\|_{\mathcal{U}} + \frac{1}{\lambda} |q|_c.
\end{align*}

For the saddle point problem with penalty, one can show the 
following inf-sup condition by means of Braess \cite{braess1996stability} 
and his dialog with Kirmse in 1990.


\begin{Proposition}[inf-sup condition for saddle point problems with penalty]\label{infsup}
Assume $g(0) = \kappa > 0$ and 
let the bilinear form be
$\mathcal{V}$-elliptic.
Then it holds the inf-sup condition
\begin{align*}
 \inf_{(u,p)\in (\mathcal{V}\times\mathcal{U})} \sup_{(w,q)\in (\mathcal{V}\times \mathcal{U})} \frac{A_{\varphi}(u,p;w,q)}{\triple(u,p)\triple \cdot \triple (w,q)\triple} \geq \beta > 0,
\end{align*}
with $\beta$ independent of $\lambda$
and $g(\varphi)$, assuming that $0\leq \frac{1}{\lambda}\leq 1$ and $\varphi\geq 1$.
\end{Proposition}

To prove this result, the following lemma is needed. The proof 
of Proposition \ref{infsup} follows after.
 \begin{lemma}\label{lemma}
  If it holds
  \begin{align}
   \frac{2 \mu a_{\varphi}(u,u)}{\|u\|_{\mathcal{U}}} +\sup_{q\in \mathcal{U}} \frac{b_{\varphi}(u,q)}{\|q\|_{\mathcal{U}} 
   + \frac{1}{\lambda}|q|_c} \geq \alpha \|u\|_{\mathcal{V}},\label{help}
  \end{align}
  or
  \begin{align}
   \sup_{(w,q) \in (\mathcal{V} \times \mathcal{W})} \frac{A_{\varphi}(u,0;w,q)}{\triple (w,q)\triple} \geq \tilde{\alpha} \|u\|_{\mathcal{V}}\quad \text{with}\ \tilde{\alpha} > 0,\label{equi}
  \end{align}
  then the inf-sup condition of Proposition \ref{infsup} follows.
   \end{lemma}
  \proofbegin{(Lemma 1)}
  First, the equivalence of (\ref{help}) and (\ref{equi}) is shown. Then, we can prove the inf-sup condition by using the two equivalent expressions.\\
  The inequality
  (\ref{help}) can be derived by (\ref{equi}) because it holds
\begin{align*}
\begin{aligned}
\alpha \|u\|_{\mathcal{V}} \leq&\ \frac{2 \mu a_{\varphi}(u,u)}{\|u\|_{\mathcal{V}}} + \sup_{q\in \mathcal{U}} \frac{b_{\varphi}(u,q)}{\|q\|_{\mathcal{U}} + \frac{1}{\lambda}|q|_c}= \frac{A_{\varphi}(u,0;u,0)}{\triple (u,0)\triple} + \sup_{q \in \mathcal{U}} \frac{A_{\varphi}(u,0;0,q)}{\triple (0,q)\triple}\\[6pt]
\leq&\ 2 \sup_{(w,q)\in(\mathcal{V}\times\mathcal{U})} \frac{A_{\varphi}(u,0;w,q)}{\triple(w,q)\triple}.
\end{aligned}
\end{align*}
Assuming $\tilde{\alpha}\geq\frac{\alpha}{2}$ implies that (\ref{help}) results from (\ref{equi}).
To prove the other direction, the
Cauchy-Schwarz inequality is applied. As a result for the definite quadratic form $a_{\varphi}(\cdot,\cdot)$, it holds $a_{\varphi}(u,w)^2\leq a_{\varphi}(u,u) \cdot a_{\varphi}(w,w)$. 
It allows the following estimate:
\begin{align*}
 \begin{aligned}
 \tilde{\alpha}\|q\|_{\mathcal{V}} \leq&\ \sup_{(w,q)\in(\mathcal{V}\times\mathcal{U})} \frac{A_{\varphi}(u,0;w,q)}{\triple(w,q)\triple}\leq \sup_{(w,q)\in(\mathcal{V}\times\mathcal{U})} \frac{2 \mu a_{\varphi}(u,w)}{\triple (w,q) \triple} + \sup_{(w,q)\in(\mathcal{V}\times\mathcal{U})} \frac{b_{\varphi}(u,q)}{\triple(w,q)\triple}\\[8pt]
=&\ \sup_{w\in \mathcal{V}} \frac{2 \mu a_{\varphi}(u,w)}{\|w\|_{\mathcal{V}}} + \sup_{q\in\mathcal{U}} \frac{b_{\varphi}(u,q)}{\triple (0,q)\triple}\leq \left[\|a\|2 \mu a_{\varphi}(u,u)\right]^{\frac{1}{2}}+ \sup_{q\in\mathcal{U}} \frac{b_{\varphi}(u,q)}{\|q\|_{\mathcal{U}} + \frac{1}{\lambda} |q|_c}\\[8pt]
  \leq&\ \frac{\|a\| 2 \mu a_{\varphi}(u,u)}{\tilde{\alpha} \|u\|_{\mathcal{V}}} + 2\sup_{q\in\mathcal{U}} \frac{b_{\varphi}(u,q)}{\|q\|_{\mathcal{U}} + \frac{1}{\lambda}|q|_c}.
 \end{aligned}
\end{align*}
With $\alpha \geq \frac{\tilde{\alpha}}{2+\frac{\|a\|}{\tilde{\alpha}}}$,
the equivalence of (\ref{help}) and (\ref{equi}) is given.\\

Next, by means of the equivalent statements, we prove that the inf-sup condition in Proposition \ref{infsup} follows:\\

Assume $(u,p)\in (\mathcal{V}\times \mathcal{U})$. For a better overview, we define
 \begin{align*}
 \text{SUP}:= \sup_{(w,q)\in(\mathcal{V}\times\mathcal{U})} \frac{A_{\varphi}(u,p;w,q)}{\triple (w,q)\triple}.
 \end{align*}
 Via the estimate
 \begin{align*}
 \begin{aligned}
A_{\varphi}(u,p;u,-p) = 2 \mu a_{\varphi}(u,u) + \frac{1}{\lambda} c_{\varphi}(p,p)\geq \frac{1}{\lambda} c_{\varphi}(p,p) = \frac{1}{\lambda}|p|_c^2,
\end{aligned}
 \end{align*}
it follows that
 \begin{align}
  \begin{aligned}
   \frac{1}{\sqrt{\lambda}}|p|_c \leq \frac{A_{\varphi}(u,p;u,-p)}{\triple(u,p)\triple} \cdot \frac{\triple(u,p)\triple}{\frac{1}{\sqrt{\lambda}} |p|_c}\leq \frac{\triple (u,p)\triple}{\frac{1}{\sqrt{\lambda}} |p|_c}\ \text{SUP}.\label{one}
  \end{aligned}
 \end{align}
The standard inf-sup condition for the saddle point problem without a penalty term provides:
\begin{align}
 \begin{aligned}
  \beta \|p\|_{\mathcal{U}}\leq \sup_{w\in \mathcal{V}} \frac{b_{\varphi}(w,p)}{\|w\|_{\mathcal{V}}}= \sup_{w\in \mathcal{V}} \frac{A_{\varphi}(u,p;w,0)-2 \mu a_{\varphi}(u,w)}{\|w\|_{\mathcal{V}}}\leq \text{SUP} + 2\mu \|a\| \|u\|_{\mathcal{V}}.\label{two}
 \end{aligned}
\end{align}
Via the estimate (\ref{equi}) it follows
\begin{align}
 \begin{aligned}
  \tilde{\alpha}\|u\|_{\mathcal{V}} \leq&\ \frac{2 \mu a_{\varphi}(u,u)}{\|u\|_{\mathcal{V}}} + \sup_{q\in \mathcal{U}} \frac{b_{\varphi}(u,q)}{\|q\|_{\mathcal{U}} + \frac{1}{\lambda}|q|_c}\leq \frac{A_{\varphi}(u,p;u,-p)}{\triple (u,p)\triple} \frac{\triple(u,p)\triple}{\|u\|_{\mathcal{V}}}\\[8pt]
  &\ + \sup_{q\in \mathcal{U}} \frac{A_{\varphi}(u,p;0,q) + \frac{1}{\lambda} c_{\varphi}(p,q)}{\triple(0,q)\triple}\leq \frac{\triple(u,p) \triple}{\|u\|_{\mathcal{V}}} \text{SUP} + \text{SUP} + \frac{1}{\sqrt{\lambda}} |p|_c.\label{three}
 \end{aligned}
\end{align}
For the next step, we use the following relation:
\begin{align}
r\leq \frac{s^2}{r} + t \Rightarrow\ r\leq s+t\quad \text{for}\ r,s,t \in \mathbb{R}^+.\label{cons}
\end{align}
Two cases have to be considered:
\begin{itemize}
 \item[$i)$] Assume that
 \begin{align*}
\frac{1}{\sqrt{\lambda}} |p|_c \leq \frac{1}{2} \tilde{\alpha} \|u\|_{\mathcal{V}}.
 \end{align*}
From (\ref{three}) and by using (\ref{two}) one can follow
\begin{align*}
\begin{aligned}
\frac{1}{2}\tilde{\alpha} \|u\|_{\mathcal{V}} \leq&\ \text{SUP} \left(\frac{\|u\|_{\mathcal{V}} + \|p\|_{\mathcal{U}} + \frac{1}{\sqrt{\lambda}} |p|_c}{\|u\|_{\mathcal{V}}} + 1 \right)\leq \text{SUP}\ \left(2 + \frac{\frac{1}{\beta} \text{SUP}}{\|u\|_{\mathcal{V}}} + \frac{ 2\mu \|a\|}{\beta} + \frac{1}{2}\tilde{\alpha} \right),\\[8pt]
\Leftrightarrow \|u\|_{\mathcal{V}} \leq&\ \text{SUP} \left(\frac{4}{\tilde{\alpha}} + \frac{\frac{2}{\tilde{\alpha}\beta} \text{SUP}}{\|u\|_{\mathcal{V}}} + \frac{4 \mu \|a\|}{\tilde{\alpha}\beta} + 1 \right).
\end{aligned}
\end{align*}
In the next step we use (\ref{cons}) with 
\begin{align*}
\begin{aligned}
r=&\ \|u\|_{\mathcal{V}},\\
s=&\ \sqrt{\frac{2}{\tilde{\alpha}\beta} \text{SUP}},\\
t=&\ \text{SUP}\ \frac{4}{\tilde{\alpha}}+\text{SUP}\ \frac{4 \mu\|a\|}{\tilde{\alpha}\beta}+\text{SUP}.
\end{aligned}
\end{align*}
Thus we obtain
\begin{align*}
\|u\|_{\mathcal{V}} \leq \text{SUP} \left(\frac{4}{\tilde{\alpha}} + \sqrt{\frac{2}{\tilde{\alpha}\beta}} + \frac{4 \mu \|a\|}{\tilde{\alpha}\beta}+1\right).
\end{align*}
Bounds for the other norms are given by (\ref{two}) and relation (\ref{cons}).
\item[$ii)$] Assume that
\begin{align}
\frac{1}{\sqrt{\lambda}} |p|_c > \frac{1}{2} \tilde{\alpha} \|u\|_{\mathcal{V}}.\label{assump}
\end{align}
First, via (\ref{two}) and(\ref{three}) it follows
\begin{align}
\|p\|_{\mathcal{U}} \leq \frac{\text{SUP}}{\beta} + \frac{2\|a\|}{\tilde{\alpha}\beta} \frac{1}{\sqrt{\lambda}} |p|_c.\label{helpeq}
\end{align}
Using (\ref{one}), in a second step using (\ref{helpeq}) and the assumption (\ref{assump}), we conclude that
\begin{align*}
\begin{aligned}
\frac{1}{\sqrt{\lambda}}|p|_c \leq&\ \text{SUP}\ \frac{\|u\|_{\mathcal{V}} + \|p\|_{\mathcal{U}}+\frac{1}{\sqrt{\lambda}} |p|_c}{\frac{1}{\sqrt{\lambda}}|p|_c}\leq \text{SUP}\ \left(\frac{2}{\tilde{\alpha}}+\frac{\frac{\text{SUP}}{\beta}}{\frac{1}{\sqrt{\lambda}}|p|_c} + \frac{2\|a\|}{\tilde{\alpha}\beta}+1\right)\\[6pt]
\leq&\ \left(1+\frac{2}{\tilde{\alpha}}+\frac{1}{\sqrt{\beta}}+\frac{2\|a\|}{\tilde{\alpha}\beta}\right)\ \text{SUP}.
\end{aligned}
\end{align*}
The bounds for the other norms can be derived by using (\ref{one}), (\ref{two}), (\ref{three}) and (\ref{cons}).
\end{itemize}
\proofend\\

\proofbegin{(Proposition 1)}
The ellipticity of $\mathcal{V}$ with 
\begin{align}
a_{\varphi}(u,u)\geq \alpha\|u\|^2_{\mathcal{V}}
\end{align}
allows to justify the validity of (\ref{help}). Via the equivalence in Lemma 1, the statement of Proposition \ref{infsup} follows and the proof is completed.
\proofend\\


Proposition \ref{infsup} establishes the stability of the problem formulation in (\ref{41}) and (\ref{42}). 

This in turn allows to formulate a stable discretization of the mixed system (\ref{mixed1}).
The discretized system reads:

Find $u_h \in \mathcal{V}_h \subset \mathcal{V}$ and $p_h \in \mathcal{U}_h\subset \mathcal{U}$ such that 
\begin{align*}
  \begin{aligned}
   2 \mu (g(\varphi)E_{\text{lin}}(u_h), E_{\text{lin}}(w)) + (g(\varphi) \nabla \cdot w, p_h) =&\ 0\quad \forall w \in \mathcal{V}_h \subset \mathcal{V},\\
(g(\varphi) \nabla\ \cdot u_h,q) - \frac{1}{\lambda} (g(\varphi) p_h,q) =&\ 0\quad \forall q \in \mathcal{U}_h \subset \mathcal{U}.
\end{aligned}
\end{align*}


\subsection{Mixed form of phase-field models}\label{newModel}
In this section, we return to the full phase-field description 
and use the previously derived mixed formulation in (\ref{41}) and (\ref{42}) for the 
displacement equation.

The complete phase-field model with a mixed form of the $u$-equation reads as follows:\\
Find $u \in \mathcal{V}$, $p \in \mathcal{U}$ and $\varphi \in \mathcal{W}$ such that 
 \begin{align*}
 \begin{aligned}
  2\mu (g(\varphi) E_{\text{lin}}(u),E_{\text{lin}}(w))+\lambda(g(\varphi)\nabla \cdot w, p)=&\ 0\quad \forall w \in \mathcal{V},\\
 (g(\varphi) \nabla \cdot u,q) - \frac{1}{\lambda} (g(\varphi) p,q)=&\ 0\quad \forall q\in \mathcal{U},
 \end{aligned}
 \end{align*}
  \begin{align*}
 \begin{aligned}
  (1- \kappa)(\varphi\ 2\mu E_{\text{lin}}(u):&\ E_{\text{lin}}(u),\psi-\varphi) + (1-\kappa)(\varphi \lambda \nabla \cdot u \textbf{I}:E_{\text{lin}}(u),\psi-\varphi)\\
  +&\ G_c (-\frac{1}{\epsilon} (1-\varphi),\psi-\varphi)+ G_c \epsilon(\nabla \varphi,\nabla(\psi-\varphi))\geq 0\quad \forall \psi \in\ \mathcal{K}\subset \mathcal{W}.
 \end{aligned}
 \end{align*}

\subsection{Stress split into tensile and compressive forces}
The next problem formulation takes into account the split of the stress tensor $\sigma(u)$ into 
tension and compression. 
The tensile stresses are named $\sigma^+(u)$, the compressive stresses summarized in $\sigma^-(u)$. The thermodynamic consistency of $\sigma^+(u)$ and $\sigma^-(u)$ has been discussed in \cite{MieWelHof10a} and \cite{PhaAmMarMau11}.
They are defined as:
\begin{align*}
\begin{aligned}
 \sigma^+(u)=&\ 2\mu E_{\text{lin}}^+(u) + \lambda \max \{0, \text{tr} (E_{\text{lin}}(u))\} \textbf{I},\\
 \sigma^-(u)=&\ 2\mu (E_{\text{lin}}(u)-E_{\text{lin}}^+(u)) + \lambda (\text{tr} (E_{\text{lin}}(u))-\max\{0,\text{tr} (E_{\text{lin}}(u))\})\textbf{I},
 \end{aligned}
\end{align*}
with $E_{\text{lin}}^+(u) := P \Lambda^+ P^T$. In 2D, the quantity $\Lambda^+$ is a two-dimensional diagonal matrix containing the eigenvalues $\lambda_1(u)$ and $\lambda_2(u)$ of the strain tensor $E_{\text{lin}}(u)$. The corresponding eigenvectors are denoted by
$v_1(u)$ and $v_2(u)$. The matrix $P$
in $E_{\text{lin}}^+(u)$ is defined as $P:= (v_1(u),v_2(u))$.\\
Beside the pressure variable $p$ derived from the mixed form with penalty, we consider tensile and compressive stresses now. For this reason, the positive part of the pressure $p^+ \in L_2(\Omega)$ has to be defined as
$p^+ := \max \{p,0\},$
such that the tensile and compressive parts of the stress tensor are reformulated to:
\begin{align*}
\begin{aligned}
 \sigma^+(u)=&\ 2\mu E_{\text{lin}}^+(u) + p^+ \textbf{I},\\
 \sigma^-(u)=&\ 2\mu (E_{\text{lin}}(u)-E_{\text{lin}}^+(u))+ (p - p^+) \textbf{I}.
 \end{aligned}
\end{align*}
Hence, the total system contains three unknown variables $u, p, \varphi$ as denoted in the following.


\begin{form}[Final mixed formulation]\label{formMixed}\hfill
Given the initial data $\varphi^0 \in \mathcal{K}$.
Find $u:= u^n \in \mathcal{V}$, $p:=p^n \in \mathcal{U}$ and $\varphi:=
\varphi^n \in \mathcal{K} \subset \mathcal{W}$ for loading steps
$n=1,2,\ldots, N$ such that 
 \begin{align*}
 \begin{aligned}
  ((1- \kappa)\varphi^2 + \kappa) (2\mu E_{\text{lin}}^+(u) + p^+ \textbf{I},\nabla w )+ (2\mu (E_{\text{lin}}(u)-E_{\text{lin}}^+(u)),E_{\text{lin}}(w))\\ + ((p - p^+)\textbf{I},E_{\text{lin}}(w)) =&\ 0\quad \forall w\in \mathcal{V},\\
 (\nabla \cdot u,q) - \frac{1}{\lambda} (p,q)=&\ 0\quad \forall q\in \mathcal{U},\\
 (1- \kappa)(\varphi 2E_{\text{lin}}^+(u) + p^+ \textbf{I} : E_{\text{lin}}(u),\psi-\varphi) + G_c (-\frac{1}{\epsilon} (1-\varphi),\psi-\varphi)\\ 
  + G_c \epsilon (\nabla \varphi,\nabla(\psi-\varphi) )\geq&\ 0\quad \forall \psi \in \mathcal{K}.
  \end{aligned}
 \end{align*}
\end{form}

Based on the formulation of the discrete problem, the
numerical steps, particularly the regularization, the discretization and the solution algorithm 
are discussed in the subsequent section.


\section{Numerical Treatment}\label{Solving}
The numerical solution proceeds from Formulation \ref{formMixed}.
Concerning robustness and efficiency, we made good experiences
treating the phase-field system in a monolithic fashion, e.g., \cite{Wi17_SISC,Wi17_CMAME}.
At first, the handling of the 
crack irreversibility constraint is clarified. In the adjacent section, the spatial
discretization and the overall solution method are explained. 


\subsection{Imposing the crack irreversibility constraint}
To realize the inequality constraint $\varphi^{n-1}\leq \varphi^n$, we introduce a Lagrange multiplier $\tau$ similar to e.g., \cite{ito2008lagrange,rockafellar1993lagrange,vexler2008adaptive}.
Using a Lagrange multiplier $\tau \in \mathcal{X}$, an additional complementarity condition
\begin{align}
\begin{aligned}
\tau\geq&\ 0\quad \text{in}\ (\Omega \times I),\\
\varphi^n - \varphi^{n-1} \leq&\ 0\quad \text{in}\ (\Omega \times I),\\
\left(\tau,\varphi^n -\varphi^{n-1}\right)=&\ 0\quad \text{in}\ (\Omega \times I),\label{tauVar}
\end{aligned}
\end{align}
has to be satisfied.
The Lagrange multiplier acts as a fourth variable in the implementation. It makes the formulation 
more expensive due to one dimension more in the system matrix.


\subsection{Spatial discretization}\label{discretization}
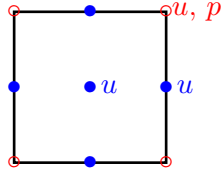
\begin{figure}[htbp!]
\centering
 \begin{tikzpicture}
  \draw (0,0) -- (2,0) -- (2,2) -- (0,2) -- cycle;
  \node[red] at (0,0) {$\circ$};
  \node[blue] at (1,0) {\textbullet};
  \node[red] at (2,0) {$\circ$};
  \node[blue] at (2,1) {\textbullet};
  \node[red] at (2,2) {$\circ$};
  \node[blue] at (1,2) {\textbullet};
  \node[red] at (0,2) {$\circ$};
  \node[blue] at (0,1) {\textbullet};
  \node[blue] at (1,1) {\textbullet};
  \node[blue] at (1.25,1) {$u$};
  \node[red] at (2.4,2) {$u$, $p$};
  \node[blue] at (2.25,1) {$u$};
 \end{tikzpicture}
\caption{Conforming quadrilateral Stokes-elements of the type $Q_2 Q_1$: $Q_2$ for the displacement variable $u$ (the filled blue and the empty red bullets) and $Q_1$ for the scalar-valued pressure variable $p$ (empty red bullets).}\label{q2q1}
\end{figure}

We employ a Galerkin finite element method for 
the spatial discretization. 
To this end, the domain $\Omega$ is partitioned into quadrilaterals.
To fulfill a discrete inf-sup condition, Taylor-Hood elements with biquadratic 
shape functions ($Q_2$) for the displacement field $u$ and bilinear shape functions ($Q_1$)
for the pressure variable $p$ are used, see Figure \ref{q2q1}.
For the definition of $Q_r, r=1,2$ elements, we refer to Ciarlet \cite{ciarlet2002finite}.

We recall (see e.g., \cite{girault2012finite}):
\begin{Proposition}[Stable Taylor-Hood elements]\label{taylor}\hfill

Taylor-Hood-elements of the type $Q_2 Q_1$ fulfill a discrete inf-sup or \selectlanguage{ngerman}Babuska"=Brezzi"=condition\selectlanguage{english} \cite{babuvska1971rate,brezzi1974existence}
 \begin{align*}
  \min_{q_h \in \mathcal{U}_h}\left\{ \max_{w_h \in \mathcal{V}_h} \frac{(q_h, g(\varphi) \nabla \cdot w_h)}{\|q_h\|\ \|\nabla w_h\|}\right\} \geq \beta_h\geq \beta >0,
 \end{align*}
with a stability constant $\beta_h$, which has to be larger than the stability constant $\beta$ on the continuous level.
\end{Proposition}
 \proofbegin
For a detailed proof we refer to Brenner and Scott \cite{brenner1994mathematical}.
 \proofend\\


For the phase-field variable $\varphi$ bilinear $Q_1$ shape functions are used. The Lagrange multiplier $\tau$ is discretized in the dual basis to the $Q_1$ space denoted by $Q_1^*$.

In our numerical tests, we need to investigate first 
the behavior with respect to higher-order finite elements. 
For this reason, we define Formulation \ref{formVariational}
using $Q_1$ elements for both $u$ and $\varphi$ as \textit{problem of reference} because 
several other groups have computed settings with this (or very similar
formulations) as well.

In our novel mixed formulation, $Q_2$ elements (biquadratic) 
for the displacement $u$ and $Q_1$ elements for the pressure variable $p$ 
are required to preserve stability.
To be more precise, 
we use $Q_2 Q_1 Q_1 Q_1^*$ elements: $Q_2$ for the 
displacement variable $u$ and $Q_1$ for 
the scalar-valued
pressure function $p$ and the phase-field function $\varphi$ and $Q_1^*$ for the Lagrange multiplier $\tau$. For the sake of a fair comparison of the reference model with the mixed model, 
the numerical results of the reference model based on Formulation \ref{formVariational} with $Q_2 Q_1$ elements are given. Numerical tests with different FE approaches are discussed in Section \ref{fem}.


\begin{form}[Discrete problem formulation]\label{formDiscrete}\hfill

Choose discrete function spaces $\mathcal{V}_h \subset \mathcal{V}$, $\mathcal{U}_h \subset \mathcal{U}$, $\mathcal{W}_h \subset \mathcal{W}$ and $\mathcal{X}_h \subset \mathcal{X}$.
Given the initial data $\varphi_h^0 \in \mathcal{W}_h$.
For the loading steps $n=1,2,\ldots, N$ solve the following system of equations:
Find $u_h \in \mathcal{V}_h$, $p_h \in \mathcal{U}_h$, $\varphi_h \in \mathcal{W}_h \subset \mathcal{W}$ and $\tau_h \in \mathcal{X}_h$ such that 
 \begin{align*}
 \begin{aligned}
   (g(\varphi)[2\mu E_{\text{lin}}^+(u_h) + p_h^+ \textbf{I}],\nabla w_h ) + (2\mu (E_{\text{lin}}(u_h)-E_{\text{lin}}^+(u_h)),E_{\text{lin}}(w_h))\\
 + ((p_h - p^+_h)\textbf{I},E_{\text{lin}}(w_h)) =&\ 0\quad \forall w_h\in \mathcal{V}_h,\\
 (g(\varphi) \nabla \cdot u_h,q_h) - \frac{1}{\lambda} (g(\varphi) p_h,q_h)=&\ 0\quad \forall q_h\in \mathcal{U}_h,\\
  (1- \kappa)(\varphi_h 2\mu E_{\text{lin}}^+(u_h) + p_h^+ \textbf{I} : E_{\text{lin}}(u_h),\psi_h) + G_c (-\frac{1}{\epsilon} (1-\varphi_h),\psi_h)\\ 
  + G_c \epsilon (\nabla \varphi_h,\nabla(\psi_h) ) + (\tau_h,\psi_h) =&\ 0\quad \forall \psi_h \in \mathcal{W}_h,\\
\tau_h\geq&\ 0,\\
\varphi_h^n - \varphi_h^{n-1} \leq&\ 0,\\
\left(\tau_h,\varphi_h^n -\varphi_h^{n-1}\right)=&\ 0.
\end{aligned}
\end{align*}
\end{form}

The last three conditions can be formulated as a semi-smooth equation utilizing a complementarity function.

\subsection{Solution algorithms}
The discrete Formulation \ref{formDiscrete} is treated in a monolithic
fashion, which renders the problem severely nonlinear (besides the
nonlinearities induced by the stress splitting and the crack irreversibility
constraint). 
To this end, we 
formulate a compact form by summing up all equations:
Given the initial data $\varphi^0$; 
for the loading steps $n=1,2,\ldots, N$:\\
Find $U_h^n:= U_h = (u_h,p_h,\varphi_h,\tau_h) \in Y_h:= (\mathcal{V}_h\times \mathcal{U}_h\times \mathcal{W}_h\times\mathcal{X}_h)$ such that
\[
A_{\varphi}(u_h,p_h,\varphi_h,\tau_h) = 0.
\]
To solve $A_{\varphi}(\cdot) = 0$, we formulate 
a residual-based Newton scheme similar to \cite{Wi17_CMAME}.
The concrete  scheme (and its implementation) 
can be found in\\ PDE/Instat/Example8 of \cite{dope,DOpElib}.
Therein, the linear system of equations is solved with a direct method provided by UMFPACK \cite{davis2004algorithm}.


\section{Numerical Tests}\label{Results}
Employing Formulation \ref{formDiscrete}, we conduct 
three studies each with two test settings in this section, which result in six numerical examples: 
\begin{itemize}
 \item Examples $1$+$2$: comparison of low-order and higher-order finite elements (Section \ref{fem})
 \item Examples $3$+$4$: mesh refinement studies (Section \ref{ref_section})
 \item Examples $5$+$6$: varying Poisson's ratio $\nu \to 0.5$ (Section \ref{poisson_section})
\end{itemize}

As first configuration, the single-edge notched shear test is considered, firstly tested with a phase-field model by Miehe et al. \cite{miehe2010phase}.
The second configuration is the L-shaped panel test proposed by Winkler \cite{winkler2001traglastuntersuchungen}, is discussed.
The programming code of both numerical tests is built on Example 8 of the instationary PDE Examples in the open-source library DOpElib \cite{dope,DOpElib}.
DOpElib in turn, is based on the deal.II finite element library \cite{dealII85}.

 
\subsection{Single edge notched pure shear test: configuration}
\begin{figure}[htbp!]
\centering
 \begin{tikzpicture}
\draw[fill=gray!30] (0,0) -- (0,5) -- (5,5) -- (5,0) -- (0,0);
\draw (0,0) -- (5,0);
\draw[red, thick] (2.5,2.5) -- (5,2.5);
\node at (3.725,2.3) {slit};
\draw (5,0) -- (5,5);
\draw[line width =2pt, blue] (5,5) -- (0,5);
\draw[-] (0,5) -- (0,0);
\draw[<->] (-0.25,0) -- (-0.25,5);
\draw[->,blue] (0.7,5.5) -- (0.1,5.5);
\node at (2.5,5.2) {$\Gamma_{\text{top}}$};
\node at (0.4,5.3) {$u_x$};
\node at (-0.9,2.5) {$10 \si{mm}$};
\draw[<->] (0,-0.5) -- (5,-0.5);
\node at (5.15,-0.5) {$x$};
\node at (-0.5,5.15) {$y$};
\node at (2.5,-0.75) {$10 \si{mm}$};

\draw[<->] (5.25,0) -- (5.25,2.5);
\node at (5.75,1.25) {$5 \si{mm}$};

\draw (0.1,0) -- (-0.05,-0.2);
\draw (0.3,0) -- (0.15,-0.2);
\draw (0.5,0) -- (0.35,-0.2);
\draw (0.7,0) -- (0.55,-0.2);
\draw (0.9,0) -- (0.75,-0.2);
\draw (1.1,0) -- (0.95,-0.2);
\draw (1.3,0) -- (1.15,-0.2);
\draw (1.5,0) -- (1.35,-0.2);
\draw (1.7,0) -- (1.55,-0.2);
\draw (1.9,0) -- (1.75,-0.2);
\draw (2.1,0) -- (1.95,-0.2);
\draw (2.3,0) -- (2.15,-0.2);
\draw (2.5,0) -- (2.35,-0.2);
\draw (2.7,0) -- (2.55,-0.2);
\draw (2.9,0) -- (2.75,-0.2);
\draw (3.1,0) -- (2.95,-0.2);
\draw (3.3,0) -- (3.15,-0.2);
\draw (3.5,0) -- (3.35,-0.2);
\draw (3.7,0) -- (3.55,-0.2);
\draw (3.9,0) -- (3.75,-0.2);
\draw (4.1,0) -- (3.95,-0.2);
\draw (4.3,0) -- (4.15,-0.2);
\draw (4.5,0) -- (4.35,-0.2);
\draw (4.7,0) -- (4.55,-0.2);
\draw (4.9,0) -- (4.75,-0.2);
 \end{tikzpicture}
 \caption{Geometry and boundary conditions of the single edge notched shear test. On the left and right side, the boundary condition in $y$-direction is $u_y = 0 \si{mm}$ and traction-free in $x$-direction. 
 On the bottom boundary it is determined $u_x = u_y = 0 \si{mm}.$ On the top boundary, it holds
 $u_y = 0 \si{mm}$ and in $x$-direction a time-dependent non-homogeneous Dirichlet condition: $u_x = t \cdot 1 \si{mm/s}.$}\label{shear_geo}
\end{figure}
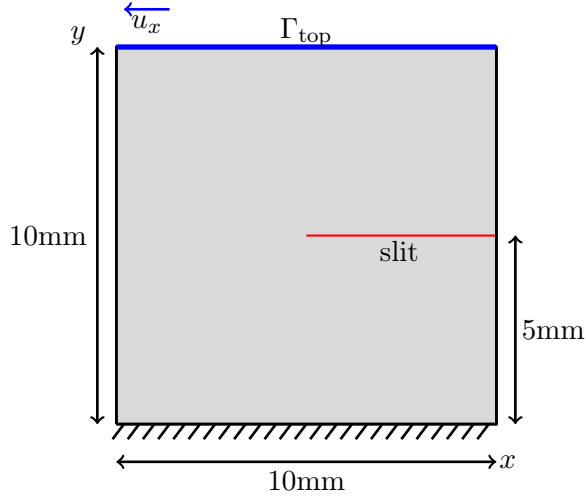

The single edge notched shear test is characterized by pure elastic crack propagation. The geometry and the material parameters are adopted from Miehe et al. \cite{miehe2010phase} and displayed in Figure \ref{shear_geo}. 
Here, the domain of interest $\Omega$ is a two-dimensional square of $10 \si{mm}$ length
with a given crack (called slit) on the right side at $5 \si{mm}$ tending to the midpoint of the square. On the bottom boundary the square is fixed, on the top boundary it is pulled with a given force. 

In Table \ref{shear_standard_parameters}, the parameters relating to material properties and parameters used for the numerical solving are listed: $\lambda_0$ and $\mu_0$ are stated such that $\nu_0=0.29999$ for the first example. 
In the following, we assume that  $\nu_0 \approx 0.3$ for the single edge notched shear test with standard settings.
The critical energy release rate $G_c$, arising in the second part of Formulation \ref{formMixed}, is defined as $G_c=2.7 \si{N/mm}$. 
The first numerical parameter in Table \ref{shear_standard_parameters} is the mesh element diameter $h_0 = 0.022 \si{mm}$. The parameter $\epsilon$ directly depends on $h_0$
with $\epsilon_0=2.0\cdot h_0$. For the standard setting we determine a loading increment of $\delta t=10^{-4} \si{s}$ and the regularization parameter $\kappa=10^{-10}$ sufficiently small.
\begin{table}[htpb!]
\centering
\renewcommand*{\arraystretch}{1.4}
\begin{tabular}{|c|c|}\hline
\multicolumn{1}{|c}{Parameter}  & \multicolumn{1}{|c|}{Value} \\ \hline \hline
   $\lambda_0$  & $121.15 \si{kN/mm^2}$\\ \hline
   $\mu_0$  & $80.77 \si{kN/mm^2}$ \\ \hline
   $\nu_0$  & $0.29999 \approx 0.3$ \\ \hline
   $G_c$  & $2.7 \si{N/mm}$ \\ \hline
    $h_0$  &  $0.022 \si{mm}$ \\ \hline
   $\epsilon_0$  & $2.0\ h_0$ \\ \hline
   $\delta t$ &  $10^{-4} \si{s}$ \\ \hline
      $\kappa$ & $10^{-10}$ \\ \hline
 \end{tabular}
\caption{Standard settings of the material and numerical parameters for the single edge notched shear test.}\label{shear_standard_parameters}
\end{table}

\subsection{L-shaped panel test: configuration}

The L-shaped panel test using a phase-field fracture model has been recently
computed by numerous groups \cite{AmGeraLoren15,bernard2012damage,feist2006embedded,meschke2007energy,unger2007modelling,Wi17_SISC}.

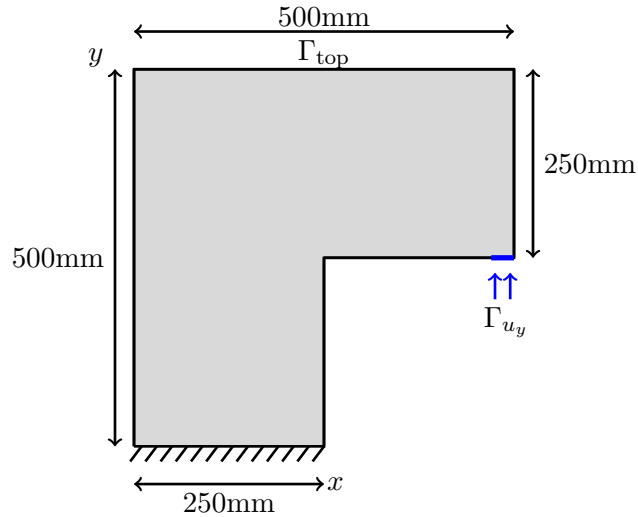
\begin{figure}[htbp!]
\centering
 \begin{tikzpicture}
\draw[fill=gray!30] (0,0) -- (2.5,0) -- (2.5,2.5) -- (5,2.5) -- (5,5) -- (0,5) -- (0,0);
\draw (2.5,0) -- (2.5,2.5);
\draw (2.5,2.5) -- (5,2.5);
\draw (5,2.5) -- (5,5);
\draw (5,5) -- (0,5);
\draw (0,5) -- (0,0);
\draw[<->] (0,5.5) -- (5,5.5);
\node at (2.5,5.7) {$500 \si{mm}$};
\node at (2.5,5.2) {$\Gamma_{\text{top}}$};
\draw[<->] (-0.25,0) -- (-0.25,5);
\node at (-0.5,5.15) {$y$};
\node at (-1,2.5) {$500 \si{mm}$};
\draw[<->] (0,-0.5) -- (2.5,-0.5);
\node at (1.25,-0.75) {$250 \si{mm}$};
\node at (2.65,-0.5) {$x$};

\draw[<->] (5.25,2.5) -- (5.25,5);
\node at (6,3.75) {$250 \si{mm}$};

\draw (0.1,0) -- (-0.05,-0.2);
\draw (0.3,0) -- (0.15,-0.2);
\draw (0.5,0) -- (0.35,-0.2);
\draw (0.7,0) -- (0.55,-0.2);
\draw (0.9,0) -- (0.75,-0.2);
\draw (1.1,0) -- (0.95,-0.2);
\draw (1.3,0) -- (1.15,-0.2);
\draw (1.5,0) -- (1.35,-0.2);
\draw (1.7,0) -- (1.55,-0.2);
\draw (1.9,0) -- (1.75,-0.2);
\draw (2.1,0) -- (1.95,-0.2);
\draw (2.3,0) -- (2.15,-0.2);
\draw (2.5,0) -- (2.35,-0.2);
\draw[line width =2pt, blue] (4.7,2.5) -- (5.0,2.5);
\draw[blue,->] (4.75,1.95) -- (4.75,2.3);
\draw[blue,->] (4.95,1.95) -- (4.95,2.3);
\node at (4.9,1.65) {$\Gamma_{u_y}$};
 \end{tikzpicture}
 \caption{Geometry and boundary conditions of the L-shaped panel test. The lower left boundary is fixed with $u_x = u_y = 0 \si{mm}.$ 
 In the right, marked corner, a special cyclic displacement condition for $u_y$ is given, defined in (\ref{uy}) and depicted in Figure \ref{CyclicLoading}.}\label{lShaped}
\end{figure}

\begin{figure}[htbp!]
\centering
\begin{tikzpicture}[xscale=0.8,yscale=0.75]
\draw[-] (-0.15,4)--(0.15,4);
\node at (-0.5,4) {0.4};
\draw[-] (-0.15,3)--(0.15,3);
\node at (-0.5,3) {0.3};
\draw[-] (-0.15,2)--(0.15,2);
\node at (-0.5,2) {0.2};
\draw[-] (-0.15,1)--(0.15,1);
\node at (-0.5,1) {0.1};
\draw[-] (-0.15,0)--(0.15,0);
\node at (-0.5,0) {0};
\draw[-] (-0.15,-1)--(0.15,-1);
\node at (-0.5,-1) {-0.1};
\draw[-] (-0.15,-2)--(0.15,-2);
\node at (-0.5,-2) {-0.2};
\draw[->] (0,0)--(7.5,0);
\draw[-] (1.5,-0.15)--(1.5,0.15);
\node at (1.5,-0.5) {300};
\draw[-] (4,-0.15)--(4,0.15);
\node at (4,-0.5) {800};
\draw[-] (7,-0.15)--(7,0.15);
\node at (6.5,-0.5) {1400};
\node at (6.7,0.35) {Load step};
\draw[->] (0,-2.2)--(0,4.25);
\node at (0,4.55) {Displacement[mm]};
\draw[blue,->] (0,0) -- (1.5,3);
\draw[blue] (0,0) -- (1.5,3) node [midway, above, sloped] (TextNode) {push};
\draw[red,dashed,->] (1.5,3) -- (4,-2);
\draw[red,dashed] (1.5,3) -- (4,-2) node [midway, above, sloped] (TextNode) {pull};
\draw[blue,->] (4,-2) -- (7,4);
\draw[blue] (4,-2) -- (7,4) node [midway, above, sloped] (TextNode) {push};
\node[blue] at (5.5,4.25) {up to $u_y=1.0$};
\draw[blue,densely dotted] (7,4) -- (7.25,4.5);
\end{tikzpicture}
\caption{The cyclic loading history on $\Gamma_{u_y}$.}\label{CyclicLoading}
\end{figure}
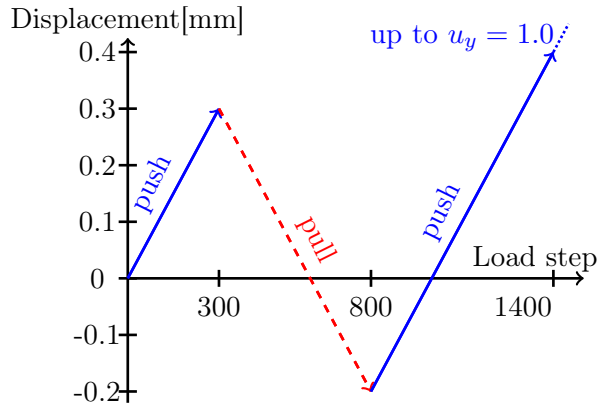

At first, the L-shaped panel test was developed by Winkler \cite{winkler2001traglastuntersuchungen} to test the crack pattern of concrete experimentally and numerically. 
Concrete is compressible with a Poisson ratio of $\nu=0.18$. To simulate fracture propagation in nearly incompressible materials, in Section \ref{poisson_section}, Poisson's ratio is increased towards the incompressible limit $\nu=0.5$.
In Figure \ref{lShaped}, the test geometry and the fitting boundary conditions of the L-shaped panel test are declared. The domain of interest has a length of $50 \si{cm}$ and resembles an $`L'$. It is fixed on the bottom part. In contrast to the first example, no initial crack is prescribed.
In the right corner $\Gamma_{u_y}$ on a small stripe of $30 \si{mm}$ at the boundary, 
a special displacement condition is defined as a loading-dependent non-homogeneous Dirichlet condition:
\begin{align}
 \begin{aligned}
 u_y = \begin{cases} t \cdot 1 \si{mm/s},\ &\text{for}\ 0.0 \si{s}\leq t < 0.3 \si{s},\\
  (0.6 - t) \cdot \si{1mm/s}, &\text{for}\ 0.3 \si{s}\leq t < 0.8 \si{s},\\
  (-1.0 + t)\cdot 1\si{mm/s}, &\text{for}\ 0.8 \si{s} \leq t < 2.0 \si{s},\label{uy}
  \end{cases}
 \end{aligned}
\end{align}
where $t$ denotes the total time. The cyclic loading defined in (\ref{uy}) is displayed in Figure \ref{CyclicLoading}.

Due to this cyclic loading the total displacement at
the end time $T = 2\si{s}$ is $1\si{mm}$.

The Lam\'{e} coefficients $\lambda_0$ and $\mu_0$ in Table \ref{l_shaped_standard_parameters} are determined such that by default $\nu_0=0.18$. With $8.9\times 10^{-5} \si{kN/mm}$
the critical energy release rate $G_c$ is determined.
The cell diameter $h_0=14.577 \si{mm}$
fixes the size of the mushy zone around the crack via $2.0\cdot h_0$ as in the
first test setup. The loading increment $\delta t$ for this test is $10^{-3}
\si{s}$ and we choose the regularization parameter $\kappa=10^{-10}$ as in the first test.
\begin{table}[htpb!]
\centering
\renewcommand*{\arraystretch}{1.4}
\begin{tabular}{|c|c|}\hline
\multicolumn{1}{|c}{Parameter}  & \multicolumn{1}{|c|}{Value} \\ \hline \hline
   $\lambda_0$  & $6.16 \si{kN/mm^2}$\\ \hline
   $\mu_0$  & $10.95 \si{kN/mm^2}$ \\ \hline
   $\nu_0$  & $0.18$ \\ \hline
   $G_c$  & $8.9\times 10^{-5} \si{kN/mm}$ \\ \hline
    $h_0$  &  $14.577 \si{mm}$ \\ \hline
   $\epsilon_0$  & $2.0\ h_0$ \\ \hline
   $\delta t$ &  $10^{-3} \si{s}$ \\ \hline
      $\kappa$ & $10^{-10}$ \\ \hline
 \end{tabular}
\caption{Standard settings of the material and numerical parameters for the L-shaped panel test.}\label{l_shaped_standard_parameters}
\end{table}


\subsection{Quantities of interest}

For both numerical tests, the functionals of interest are the \selectlanguage{ngerman}load"=displacement\selectlanguage{english} curves and the crack
path by observing the behavior of the phase-field function. In addition, we plot 
the phase-field variable at certain time steps similar to \cite{AmGeraLoren15,Wi17_SISC}.
Relative to the \selectlanguage{ngerman}load"=displacement\selectlanguage{english} curves, the load vector on the top boundary is evaluated via 
\begin{align}
 (F_x,F_y):=\int_{\Gamma_{\text{top}}} \sigma(u) n\ ds,\label{loading}
\end{align}
with the stress tensor $\sigma(u) := 2 E_{\text{lin}}(u) + \lambda \text{tr}  (E_{\text{lin}}(u)) \textbf{I}$ and the normal vector $n$. 
In the \selectlanguage{ngerman}load"=displacement\selectlanguage{english} 
curves the loading is displayed versus the displacements, which vary over time.
Within the single edge notched shear test we are particularly interested in the loading force $F_x$, 
in the frame of the L-shaped panel test we are interested in the loading force $F_y$ on $\Gamma_{\text{top}}$.
Using the definition in (\ref{loading}), the \selectlanguage{ngerman}load"=displacement\selectlanguage{english} curves of all executed numerical tests are plotted in the following sections. 
Section \ref{fem} opens a discussion on the choice and
influence of different finite element approaches. The next section contains numerical results of both tests considering meshes of different size. 
Finally in Section \ref{poisson_section}, further numerical results are presented concerning different Poisson ratios.


\subsection{The choice of finite elements}\label{fem}
To fulfill the discrete inf-sup condition, we work with the Taylor-Hood
element; see Proposition \ref{taylor}. However, higher-order finite elements, namely
$Q_2$ have not been tested in detail in the published literature.
For this reason, we first use the classical Formulation \ref{formVariational} 
and employ different finite element combinations. 

The studies with different FE approaches are performed with standard settings 
provided in the Tables \ref{shear_standard_parameters} and \ref{l_shaped_standard_parameters}.


\subsubsection{Results of the single edge notched shear test}

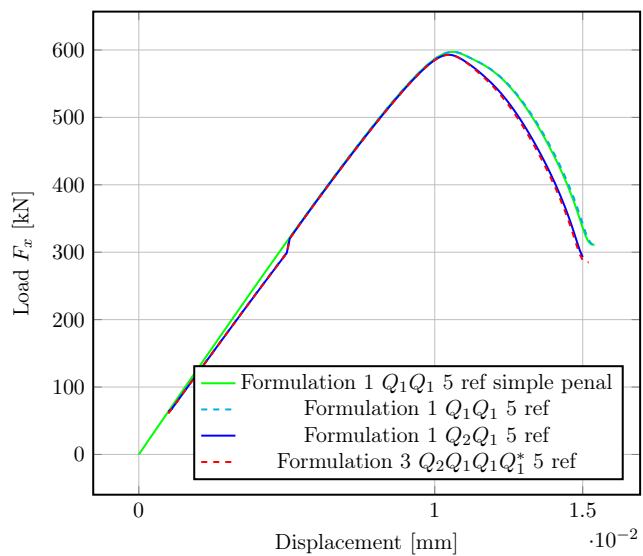
\begin{figure}[htbp!]
\centering
\begin{tikzpicture}[xscale=0.75,yscale=0.75]
\begin{axis}[
    ylabel = Load $F_x$ $\lbrack\si{kN}\rbrack$,
    xlabel = Displacement $\lbrack\si{mm}\rbrack$,
 legend pos=south east, grid =major,
    x post scale = 1.4,
    y post scale = 1.5,
  xtick={0,0.01,0.015,0.02,0.03},  
  ytick={-100,0,100,200,300,400,500,600,700}
  ]
\addplot[green]
table[x=Time,y=q1q1_5global_simplePenal,col sep=comma] {Example8_q1q1_obstacle_miehe_shear_stress.csv}; 
\addlegendentry{Formulation 1 $Q_1 Q_1$ 5 ref simple penal}
\addplot[cyan,dashed]
table[x=Time,y=q1q1_5global,col sep=comma] {Example8_q1q1_obstacle_miehe_shear_stress.csv}; 
\addlegendentry{Formulation 1 $Q_1 Q_1$ 5 ref}
\addplot[blue]  
table[x=Time,y=q2q1_5global,col sep=comma] {Example8_q2q1_obstacle_miehe_shear_stress.csv};
\addlegendentry{Formulation 1 $Q_2 Q_1$ 5 ref}
\addplot[red,dashed]  
table[x=Time,y=030_5global,col sep=comma] {030_load_displacement_Miehe.csv};
\addlegendentry{Formulation 3 $Q_2 Q_1 Q_1 Q_1^*$ 5 ref}
\end{axis}
\end{tikzpicture}
\caption{Load-displacement curves for the single edge notched shear test with $5$ steps of uniform refinement for the original implementation with $Q_2 Q_1$ and $Q_1 Q_1$ elements in comparison to the new model $Q_2 Q_1 Q_1 Q_1^*$.}\label{fem_load_displ_shear}
\end{figure}

In Figure \ref{fem_load_displ_shear} three curves are depicted: the first corresponds to the results of the reference model with $Q_1 Q_1$ elements for $u$ and the phase-field function $\varphi$, the second curve results of the same phase-field model but with $Q_2 Q_1$ elements.
The third \selectlanguage{ngerman}load"=displacement\selectlanguage{english} curve provides the data of the new model with Taylor-Hood stable $Q_2 Q_1$ elements for the displacements $u$ and the pressure variable $p$, respectively. 
The \selectlanguage{ngerman}load"=displacement\selectlanguage{english} curves with $Q_1 Q_1$ or $Q_2 Q_1$ elements or the new implementation of the mixed form with $Q_2 Q_1$ elements have a very similar course. 
In particular, results from the phase-field model with $Q_2 Q_1$ elements and the mixed phase-field model based on Formulation \ref{formDiscrete} with $Q_2 Q_1 Q_1 Q_1^*$ elements are presented.

Note that the load displacement curves, using the model implemented in Example $8$ in the folder Examples/PDE/InstatPDE of the DOpElib library with $Q_1 Q_1$ elements, are conform to the plots given in the literature, e.g., \cite{AmGeraLoren15} and \cite{Wi17_SISC}. 
The loading force is increasing until a point, where the critical energy release rate is reached and the material cracks. The crack propagation releases energy and establishes a decreasing loading after a certain maximal point until the material is broken to the bottom left corner.


\subsubsection{Results of the L-shaped panel test}
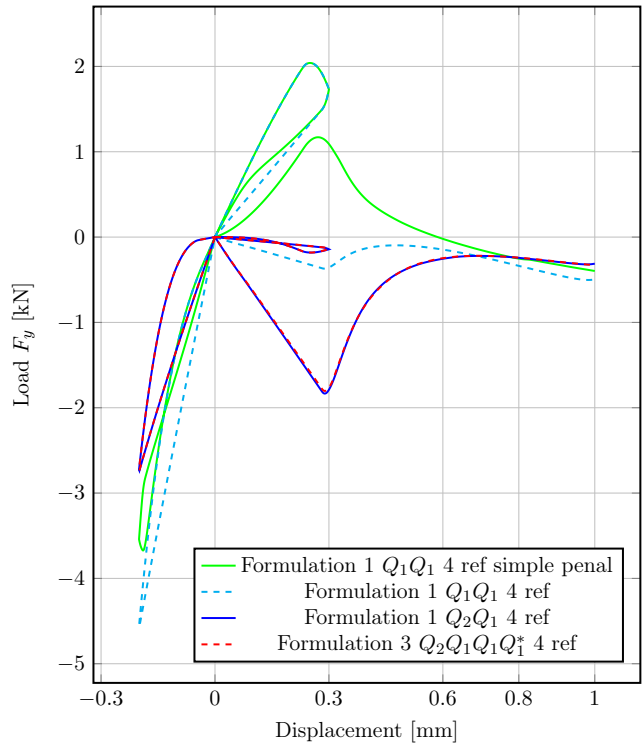
\begin{figure}[htbp!]
\centering
\begin{tikzpicture}[xscale=0.75,yscale=0.75]
\begin{axis}[
    ylabel = Load $F_y$ $\lbrack\si{kN}\rbrack$,
    xlabel = Displacement $\lbrack\si{mm}\rbrack$,
 legend pos=south east, grid =major,
    x post scale = 1.4,
    y post scale = 2.1,
  xtick={-0.3,0,0.3,0.6,0.8,1.0}, 
  ytick={-8,-7,-6,-5,-4,-3,-2,-1,0,1,2,3,4}
  ]
\addplot[green]
table[x=Loading,y=q1q1_4global_simplePenal,col sep=comma] {Example8_q1q1_load_displacement_l_shaped.csv}; 
\addlegendentry{Formulation 1 $Q_1 Q_1$ 4 ref simple penal}
\addplot[cyan,dashed]
table[x=Loading,y=q1q1_4global,col sep=comma] {Example8_q1q1_load_displacement_l_shaped.csv}; 
\addlegendentry{Formulation 1 $Q_1 Q_1$ 4 ref}
\addplot[blue]  
table[x=Loading,y=q2q1_4global,col sep=comma] {Example8_q2q1_load_displacement_l_shaped.csv};
\addlegendentry{Formulation 1 $Q_2 Q_1$ 4 ref}
\addplot[red, dashed]  
table[x=Loading,y=018_4global,col sep=comma] {018_load_displacement_l_shaped.csv};
\addlegendentry{Formulation 3 $Q_2 Q_1 Q_1 Q_1^*$ 4 ref}
\end{axis}
\end{tikzpicture}
\caption{Load-displacement curves of the L-shaped panel test for $\nu=0.18$ with $4$ steps of uniform refinement for the original implementation (Formulation \ref{formVariational}) with $Q_1 Q_1$ elements and a simple penalization method compared to $Q_1 Q_1$ elements and the new penalization strategy (Lagrange multiplier). 
Further, the load-displacement curves from the phase-field model with $Q_2 Q_1$ elements and the mixed phase-field model (Formulation \ref{formDiscrete}) with $Q_2 Q_1 Q_1 Q_1^*$ elements are presented to see the influence of different finite element approaches.}\label{fem_load_displ_l_shaped}
\end{figure}

We continue with the L-shaped panel test, which differs in a cyclic loading function and no existing crack at the beginning of the simulation.
In Figure \ref{fem_load_displ_l_shaped}, the \selectlanguage{ngerman}load"=displacement\selectlanguage{english} curves for the L-shaped panel test with different finite elements are displayed. 
The \selectlanguage{ngerman}load"=displacement\selectlanguage{english} curve resulting from $Q_1 Q_1$ elements with a simple penalization based on Formulation \ref{formVariational} agrees to the results presented in \cite{AmGeraLoren15} 
and \cite{Wi17_SISC}. 
But, as it can be seen in the curves where $Q_2$ elements are used for the
displacement field $u$, the
\selectlanguage{ngerman}load"=displacement\selectlanguage{english} values all
are negative and differ significantly from the $Q_1$ discretization.
Furthermore, the usage of a different penalization method, changes the
course of the curve, especially in the second period of pressing on the small boundary $\Gamma_{u_y}$.

\subsection{Spatial mesh refinement}\label{ref_section}
Here, mesh refinement studies are performed with a varying mesh size parameter $h$.
The finite element approximation of the following tests is based on Formulation \ref{formDiscrete} with $Q_2 Q_1 Q_1 Q_1^*$ elements.


\subsubsection{Results of the single edge notched shear test}
Beginning with different meshes for the single edge \selectlanguage{ngerman}notched\selectlanguage{english} shear test, Table \ref{shear_refine_table} presents three test cases with 4,5 and 6 steps of uniform refinement, the corresponding number of degrees of freedom (dofs) and
the bandwidth $\epsilon$ of the transition zone dependent on the maximal cell length $h$.

\begin{table}[htbp!]
\centering
\renewcommand*{\arraystretch}{1.4}
\begin{tabular}{|c|r|c|c|}\hline
\multicolumn{1}{|c}{$\#$ref} & \multicolumn{1}{|c}{$\#$dofs} & \multicolumn{1}{|c}{$\epsilon$} & \multicolumn{1}{|c|}{$h$} \\ \hline \hline
4 & 1024  & $0.088 \si{mm}$ & $0.044 \si{mm}$\\ \hline
5 & 12771  & $0.044 \si{mm}$ & $0.022 \si{mm}$ \\ \hline
6 & 50115 & $0.022 \si{mm}$ & $0.011 \si{mm}$ \\ \hline
 \end{tabular}
\caption{Values of $\epsilon$ and the mesh element diameter $h$ for 4,5 and 6 steps of uniformly refined meshes provided for the single edge notched shear test.}\label{shear_refine_table}
\end{table}

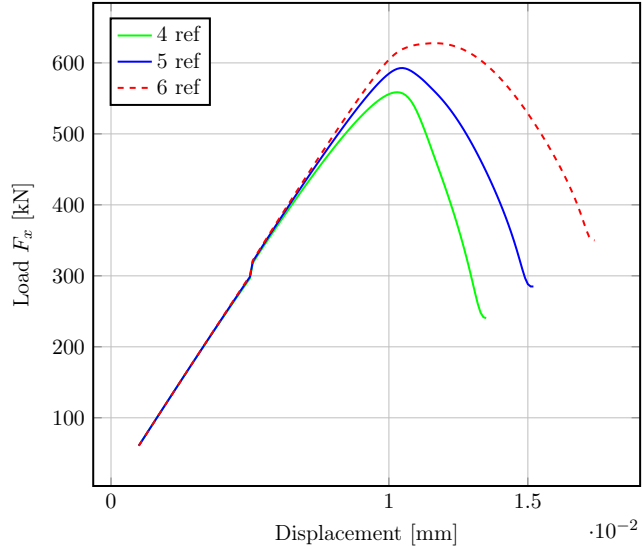
\begin{figure}[htbp!]
\centering
\begin{tikzpicture}[xscale=0.75,yscale=0.75]
\begin{axis}[
    ylabel = Load $F_x$ $\lbrack\si{kN}\rbrack$,
    xlabel = Displacement $\lbrack\si{mm}\rbrack$,
 legend pos=north west, grid =major,
    x post scale = 1.4,
    y post scale = 1.5,
  xtick={0,0.01,0.015,0.02,0.03}, 
  ytick={-100,0,100,200,300,400,500,600,700}
  ]
\addplot[green]
table[x=Time,y=030_4global,col sep=comma] {030_load_displacement_Miehe.csv}; 
\addlegendentry{4 ref}
\addplot[blue]  
table[x=Time,y=030_5global,col sep=comma] {030_load_displacement_Miehe.csv};
\addlegendentry{5 ref}
\addplot[red,dashed]  
table[x=Time,y=030_6global,col sep=comma] {030_load_displacement_Miehe.csv};
\addlegendentry{6 ref}
\end{axis}
\end{tikzpicture}
\caption{Load-displacement curves for the single edge notched shear test with 4,5 and 6 steps of uniform refinement. Poisson's ratio $\nu = 0.3$.}\label{shear_refine_plot}
\end{figure}

In Figure \ref{shear_refine_plot}, one can observe the course of the \selectlanguage{ngerman}load"=displacement\selectlanguage{english} curves with different levels of global refinement.


\subsubsection{Results of the L-shaped panel test}
The boundary condition defined in (\ref{uy}) and displayed in Figure \ref{CyclicLoading}, seems to be responsible for the course of the \selectlanguage{ngerman}load"=displacement\selectlanguage{english} curves. 
Table \ref{l_shaped_refine_table} provides all used widths $\epsilon$ depending on the mesh cell diameter $h$ for the L-shaped panel test.

\begin{table}[htbp!]
\centering
\renewcommand*{\arraystretch}{1.4}
\begin{tabular}{|c|r|r|r|}\hline
\multicolumn{1}{|c}{$\#$ref} & \multicolumn{1}{|c}{$\#$dofs} & \multicolumn{1}{|c}{$\epsilon$} & \multicolumn{1}{|c|}{$h$} \\ \hline \hline
2 & 1200  & $29.154 \si{mm}$ & $14.577 \si{mm}$\\ \hline
3 & 4800  & $14.577 \si{mm}$ & $7.289 \si{mm}$ \\ \hline
4 & 19200 & $7.289 \si{mm}$ & $3.644 \si{mm}$ \\ \hline
5 & 76800 & $3.644 \si{mm}$ & $1.822 \si{mm}$ \\ \hline
 \end{tabular}
\caption{Values of $\epsilon$ and the mesh element diameter $h$ for $2,3,4$ and $5$ steps of uniformly refined meshes provided for the L-shaped panel test.}\label{l_shaped_refine_table}
\end{table}

\begin{figure}[htbp!]
\centering
\begin{tikzpicture}[xscale=0.75,yscale=0.75]
\begin{axis}[
    ylabel = Load $F_y$ $\lbrack\si{kN}\rbrack$,
    xlabel = Displacement $\lbrack\si{mm}\rbrack$,
 legend pos=south east, grid =major,
    x post scale = 1.4,
    y post scale = 2.3,
  xtick={-0.3,0,0.3,0.6,0.8,1.0}, 
  ytick={-8,-7,-6,-5,-4,-3,-2,-1,0}
  ]
\addplot[green]
table[x=Loading,y=018_2global,col sep=comma] {018_load_displacement_l_shaped.csv}; 
\addlegendentry{2 ref}
\addplot[cyan,dashed]
table[x=Loading,y=018_3global,col sep=comma] {018_load_displacement_l_shaped.csv}; 
\addlegendentry{3 ref}
\addplot[blue]  
table[x=Loading,y=018_4global,col sep=comma] {018_load_displacement_l_shaped.csv};
\addlegendentry{4 ref}
\addplot[red,dashed]  
table[x=Loading,y=018_5global,col sep=comma] {018_load_displacement_l_shaped.csv};
\addlegendentry{5 ref}
\end{axis}
\end{tikzpicture}
\caption{Load-displacement curves of the L-shaped panel test for $\nu=0.18$ with $2,3,4$ and $5$ steps of uniform refinement.}\label{l_shaped_refine_plot}
\end{figure}
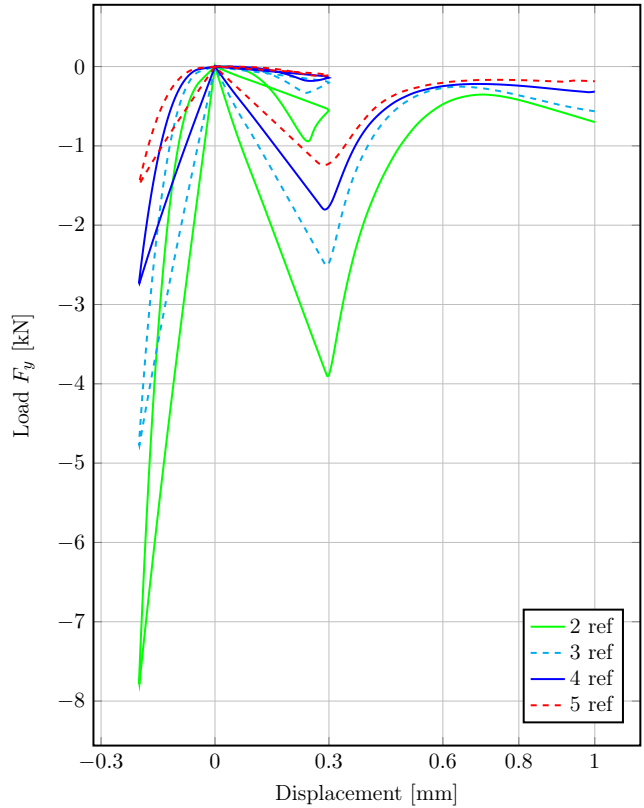

Figure \ref{l_shaped_refine_plot} provides the \selectlanguage{ngerman}load"=displacement\selectlanguage{english} curves,
corresponding to meshes with different levels of uniform refinement. 
The plotted curves significantly differ from the \selectlanguage{ngerman}load"=displacement\selectlanguage{english} curves, which can be found e.g., in \cite{AmGeraLoren15}. 
We justified in Section \ref{fem} that the reason is the choice of the FE combination.
The huge sensitivity of the \selectlanguage{ngerman}load"=displacement\selectlanguage{english} curves with respect to $h$ was also observed in \cite{Wi17_SISC} (see in particular reference $[26]$ therein).

\subsection{Variations in Poisson's ratio}\label{poisson_section}
We now turn our attention to the key objective of this work and
test the new phase-field model with Poisson ratio's towards $\nu = 0.5$.
The relevant \selectlanguage{ngerman}load"=displacement\selectlanguage{english} curves are given in the following for both introduced tests and Poisson ratios from $0.18$ to $0.4999$. 
We discuss the numerical results via
\selectlanguage{ngerman}load"=displacement\selectlanguage{english} curves and
via illustration of the crack path of the phase-field function at certain time steps with higher values of the Lam\'{e} coefficient $\lambda$. 


\subsubsection{Results of the single edge notched shear test}
In Table \ref{shear_nu_table}, 
different $\nu$-values and the corresponding values for the Lam\'{e} coefficients are listed.
A Poisson ratio of $\nu=0.3$ corresponds to the standard setting of the single edge notched shear test. 

\begin{table}[htbp!]
\centering
\renewcommand*{\arraystretch}{1.4}
\begin{tabular}{|l|c|r|}\hline
\multicolumn{1}{|c}{$\nu$} & \multicolumn{1}{|c}{$\mu$} & \multicolumn{1}{|c|}{$\lambda$} \\ \hline \hline
$0.3$  & $80.77\cdot 10^3$ & $121.15\cdot 10^3$ \\ \hline
$0.45$ & $80.77\cdot 10^3$ & $726.93\cdot 10^3$ \\ \hline
$0.49$ & $80.77\cdot 10^3$ & $3957.73\cdot 10^3$ \\ \hline
$0.499$ & $80.77\cdot 10^3$ & $40304.20\cdot 10^3$ \\ \hline
$0.4999$ & $80.77\cdot 10^3$ & $403769.00\cdot 10^3$ \\ \hline
 \end{tabular}
\caption{Tests with different Poisson's ratios approximating $\nu=0.5$ for the single edge notched shear test.}\label{shear_nu_table}
\end{table}

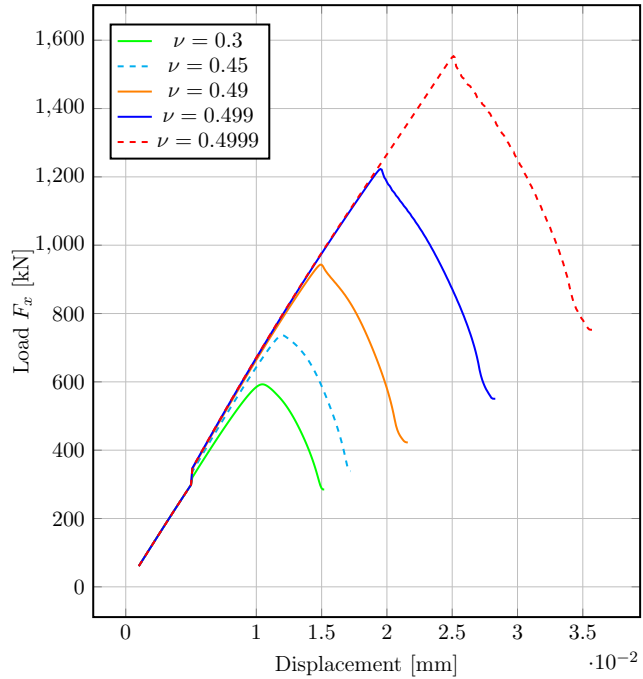
\begin{figure}[htbp!]
\centering
\begin{tikzpicture}[xscale=0.75,yscale=0.75]
\begin{axis}[
    ylabel = Load $F_x$ $\lbrack\si{kN}\rbrack$,
    xlabel = Displacement $\lbrack\si{mm}\rbrack$,
 legend pos=north west, grid =major,
    x post scale = 1.4,
    y post scale = 1.9,
  xtick={0,0.01,0.015,0.02,0.025,0.03,0.035,0.04}, 
  ytick={0,200,400,600,800,1000,1200,1400,1600}
  ]
\addplot[green]
table[x=Time,y=030_5global,col sep=comma] {nuplot_load_displacement_Miehe.csv}; 
\addlegendentry{$\nu = 0.3$}
\addplot[cyan,dashed]  
table[x=Time,y=045_5global,col sep=comma] {nuplot_load_displacement_Miehe.csv};
\addlegendentry{$\nu = 0.45$}
\addplot[orange]  
table[x=Time,y=049_5global,col sep=comma] {nuplot_load_displacement_Miehe.csv};
\addlegendentry{$\nu = 0.49$}
\addplot[blue]  
table[x=Time,y=0499_5global,col sep=comma] {nuplot_load_displacement_Miehe.csv};
\addlegendentry{$\nu = 0.499$}
\addplot[red,dashed]  
table[x=Time,y=04999_5global,col sep=comma] {nuplot_load_displacement_Miehe.csv};
\addlegendentry{$\nu = 0.4999$}
\end{axis}
\end{tikzpicture}
\caption{Load-displacement curves for the single edge notched shear test with different Poisson ratios and $5$ steps of uniform refinement.}\label{shear_nu_plot}
\end{figure}

Figure \ref{shear_nu_plot} displays the \selectlanguage{ngerman}load"=displacement\selectlanguage{english} curves with different values of the Poisson ratio $\nu$. 
With an increasing $\nu$, the loading values seem to be higher in general. The curves have a sharper maximal loading and the crack progresses later in time.
The more incompressible a material is, the more robust it appears to be against displacement forces.

\begin{figure}[htbp!]
\begin{minipage}{0.48\textwidth}
 \includegraphics[width=0.3\textwidth]{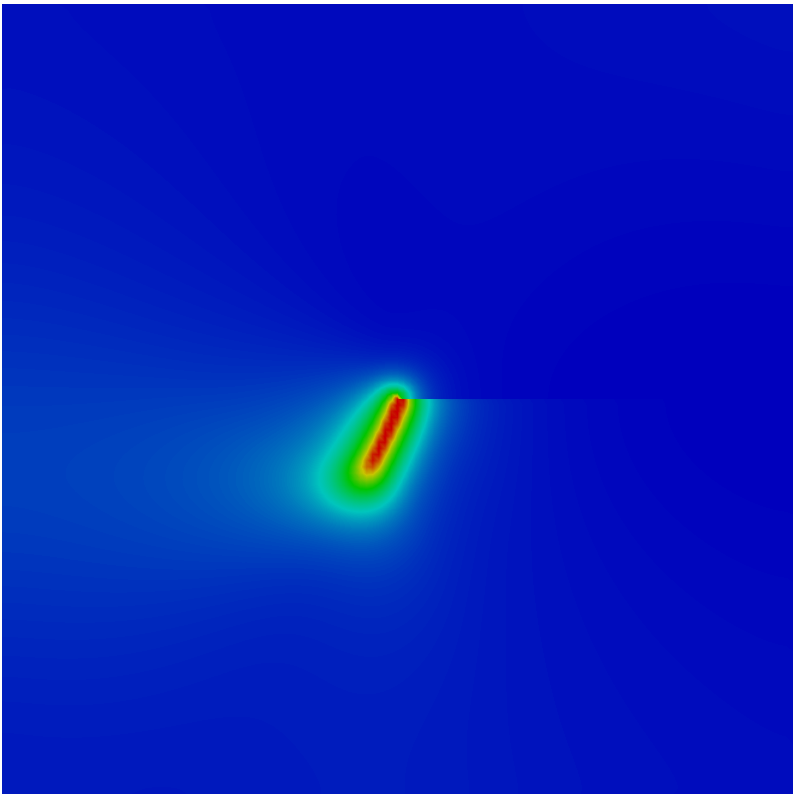}
 \includegraphics[width=0.3\textwidth]{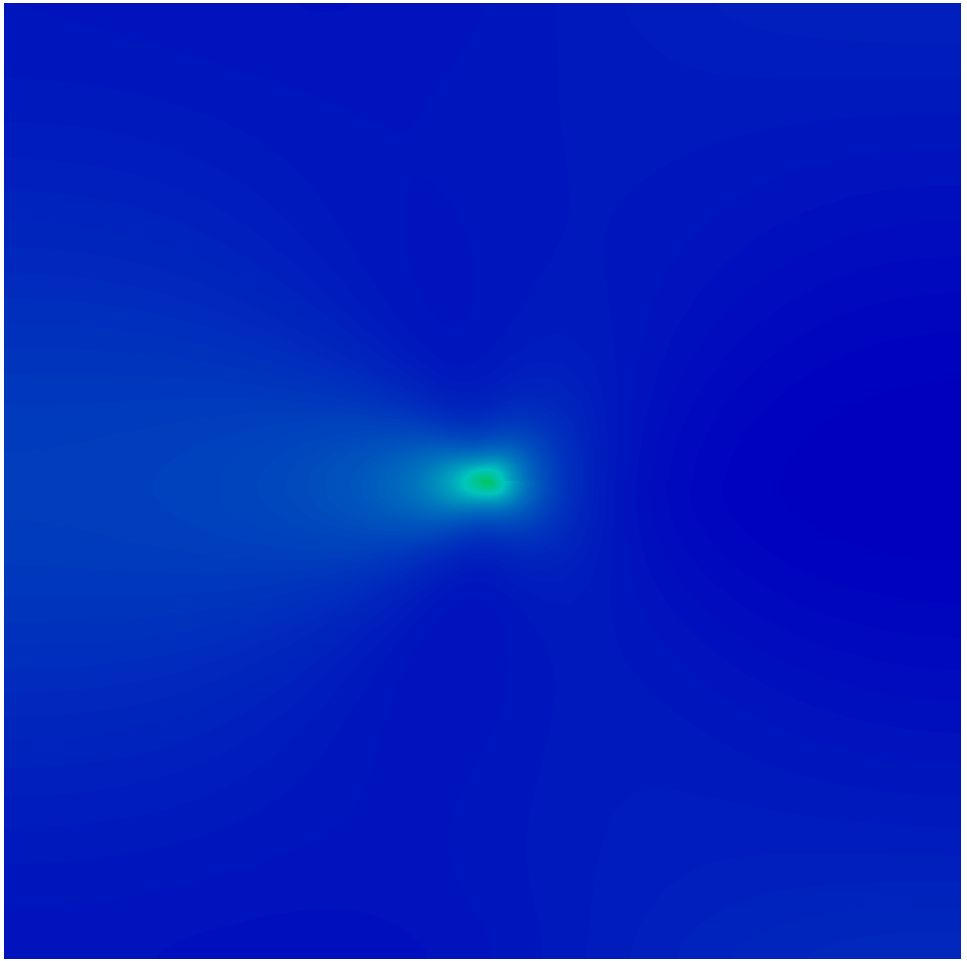}
 \includegraphics[width=0.378\textwidth]{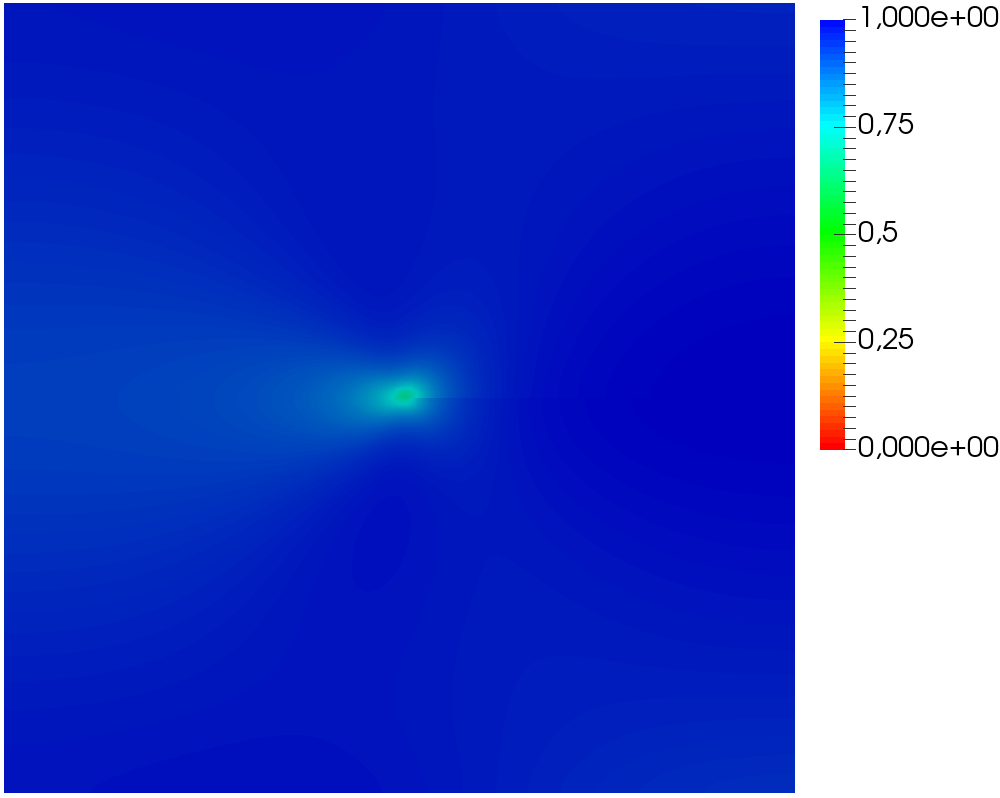}
\end{minipage}
\hfill
\begin{minipage}{0.48\textwidth}
 \includegraphics[width=0.3\textwidth]{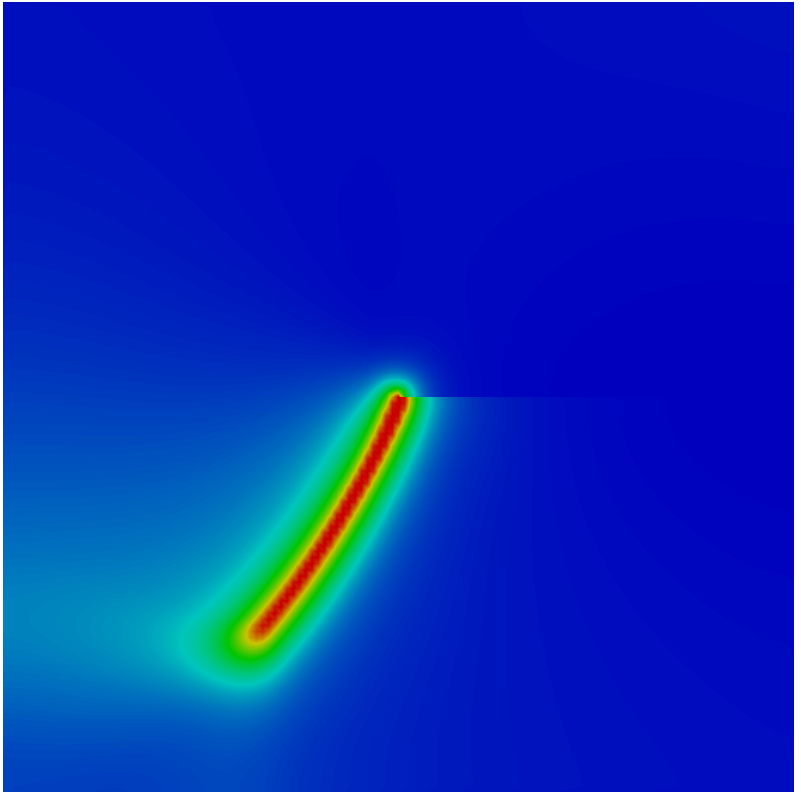}
 \includegraphics[width=0.3\textwidth]{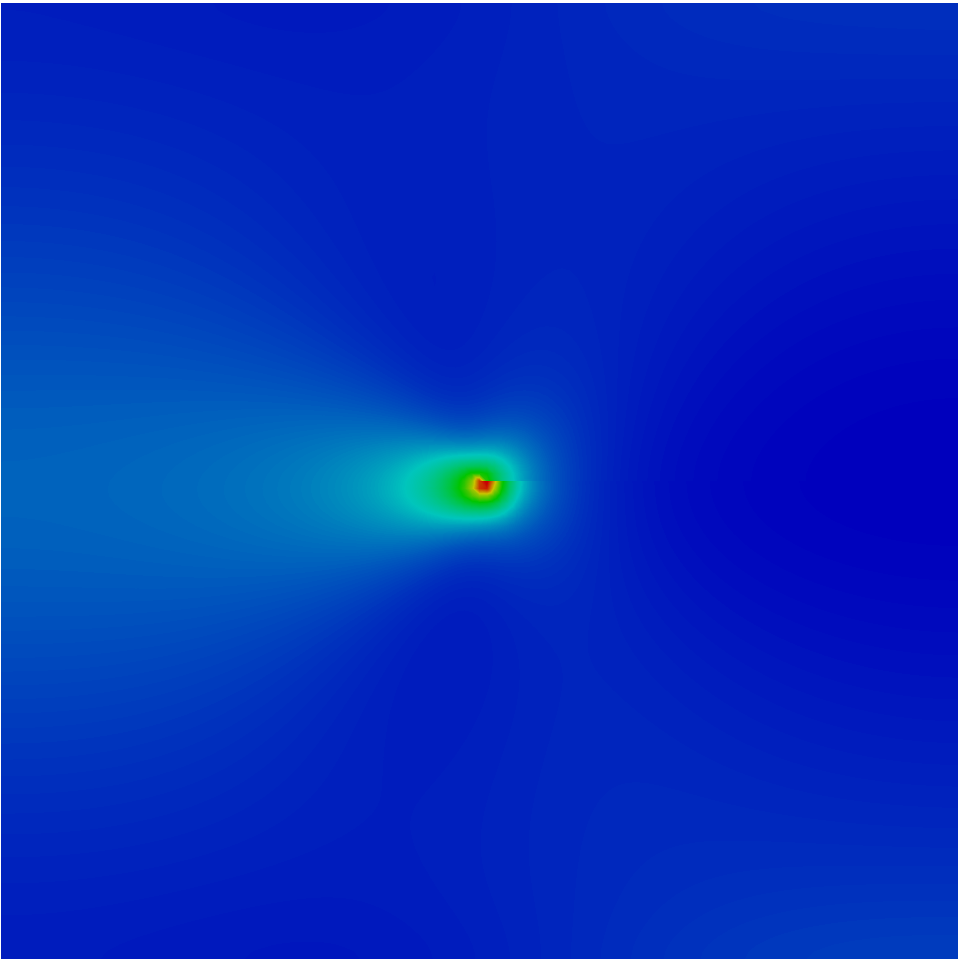}
 \includegraphics[width=0.378\textwidth]{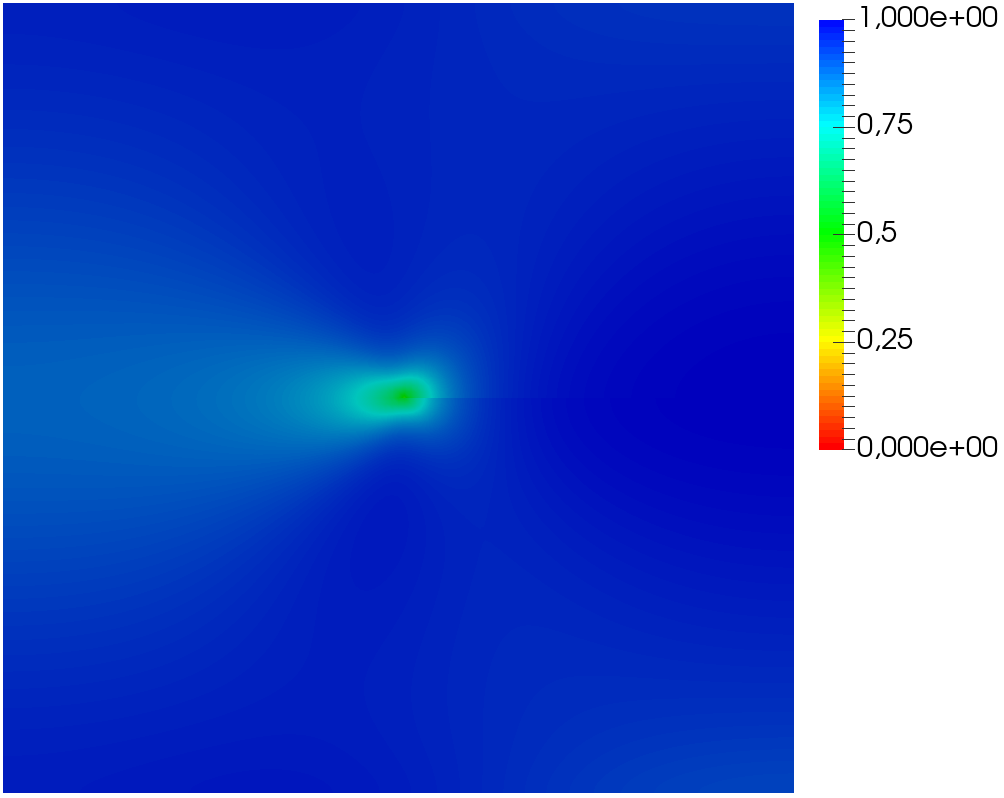}
 \end{minipage}
 \hfill
 \begin{minipage}{0.48\textwidth}
 \includegraphics[width=0.3\textwidth]{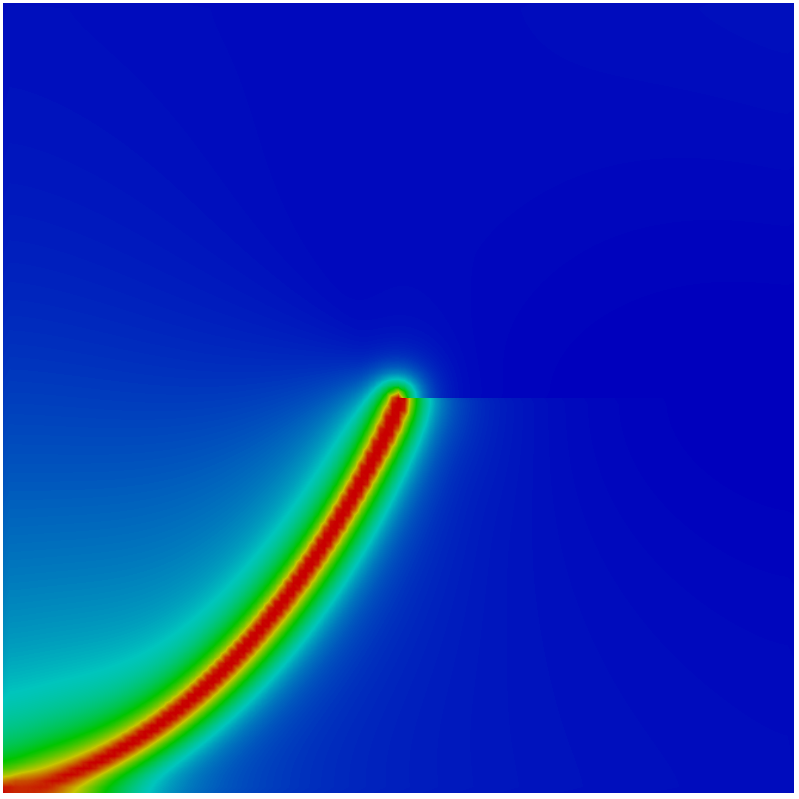}
 \includegraphics[width=0.3\textwidth]{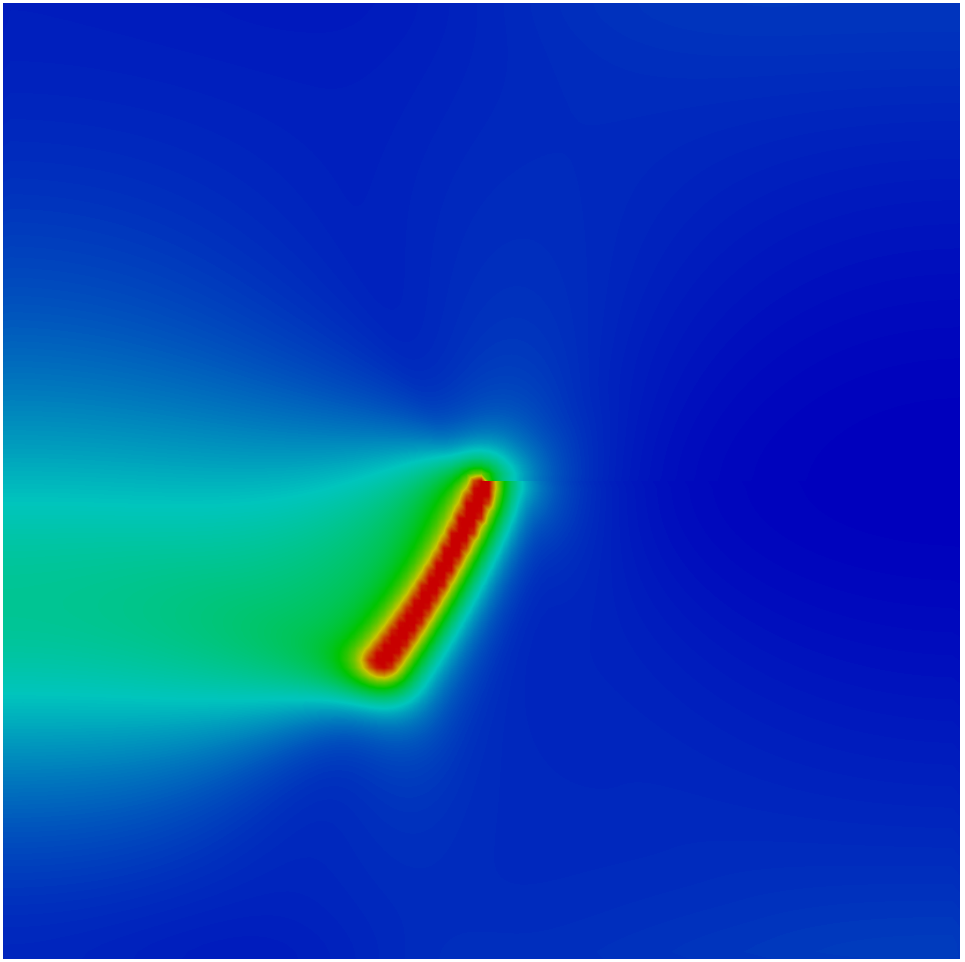}
 \includegraphics[width=0.378\textwidth]{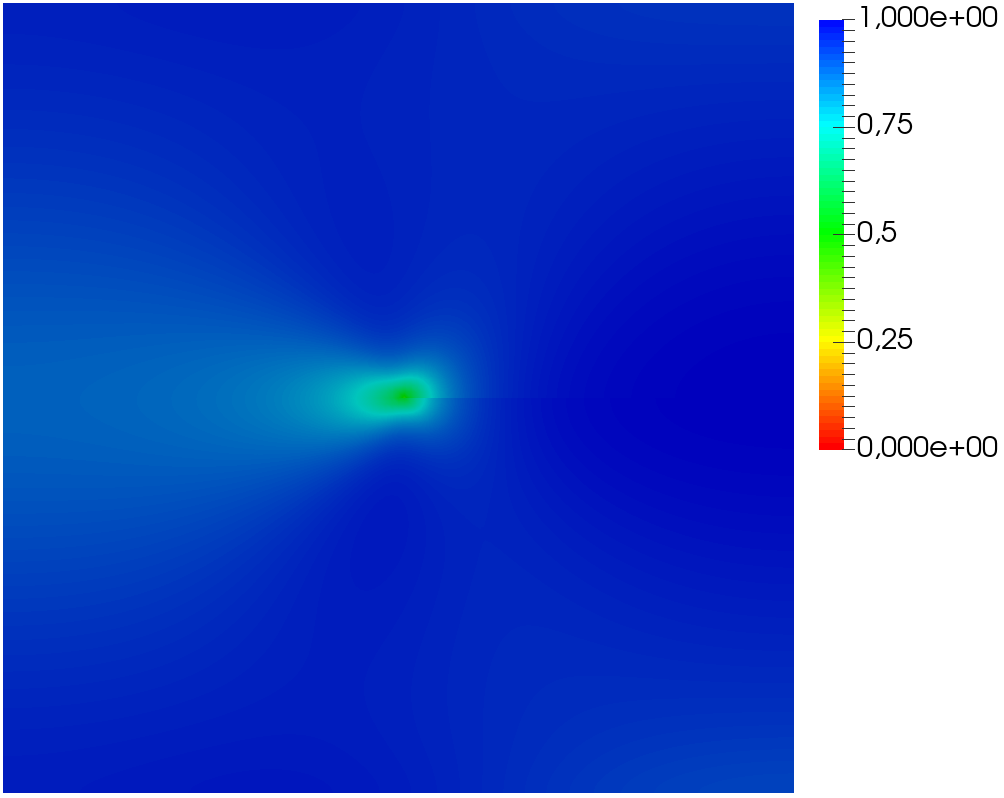}
\end{minipage}
\hfill
\begin{minipage}{0.48\textwidth}
 \includegraphics[width=0.3\textwidth]{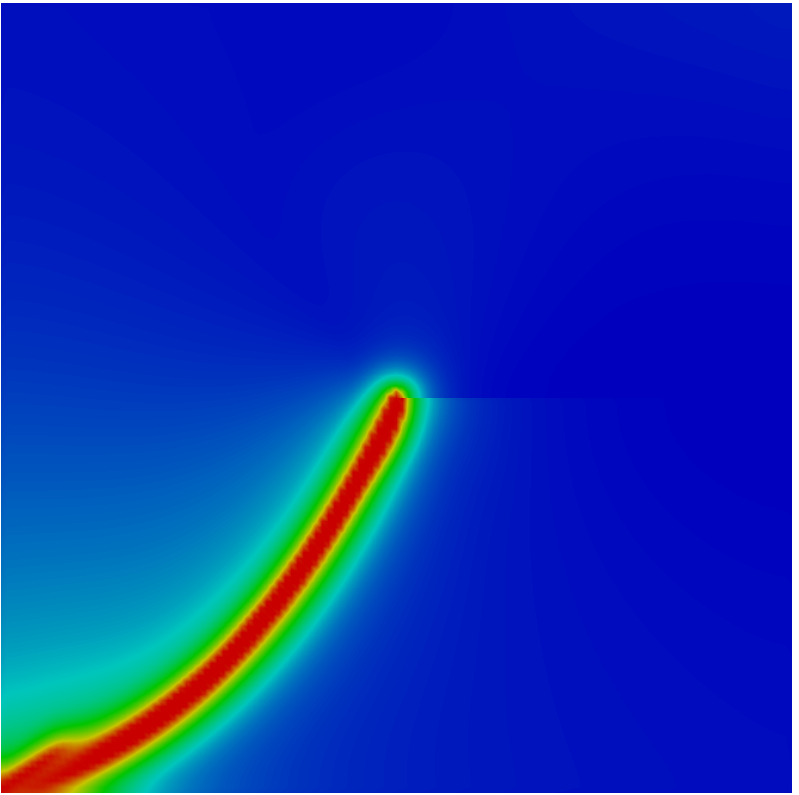}
 \includegraphics[width=0.3\textwidth]{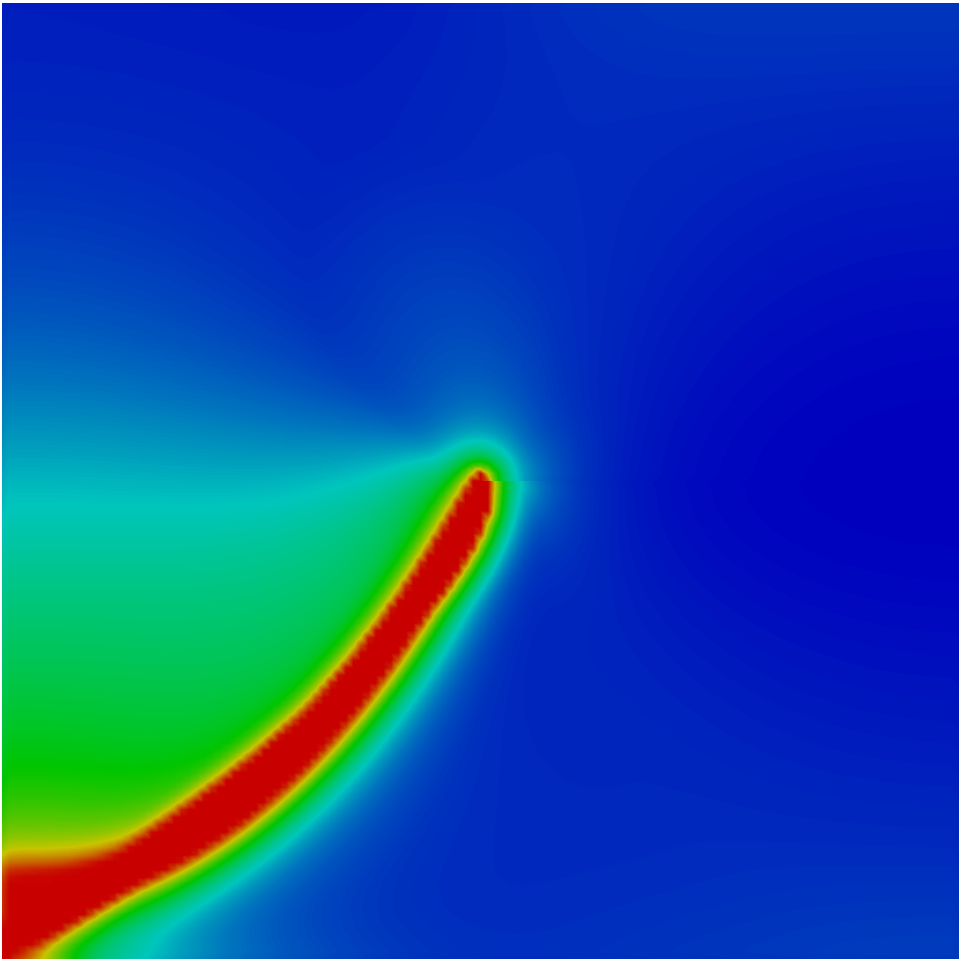}
 \includegraphics[width=0.378\textwidth]{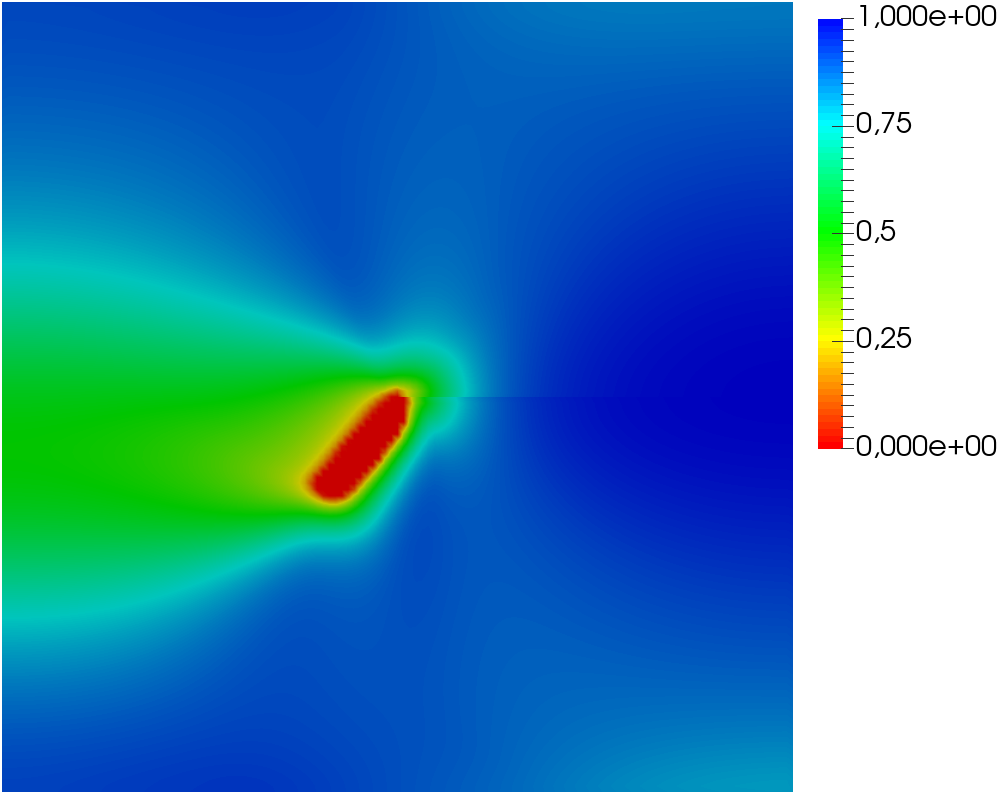}
 \end{minipage}
 \caption{Phase-field function with loading $u_y = 0.012, 0.015, 0.02$ and $0.03 \si{mm}$ from top to bottom line and for $\nu=0.30$ (left), for $\nu=0.49$ (middle) and for $\nu=0.4999$ (right) with $6$ steps of uniform refinement. }\label{shear_nu_screenshots}
\end{figure}

\begin{figure}[htbp!]
\begin{minipage}{0.48\textwidth}
 \includegraphics[width=0.432\textwidth]{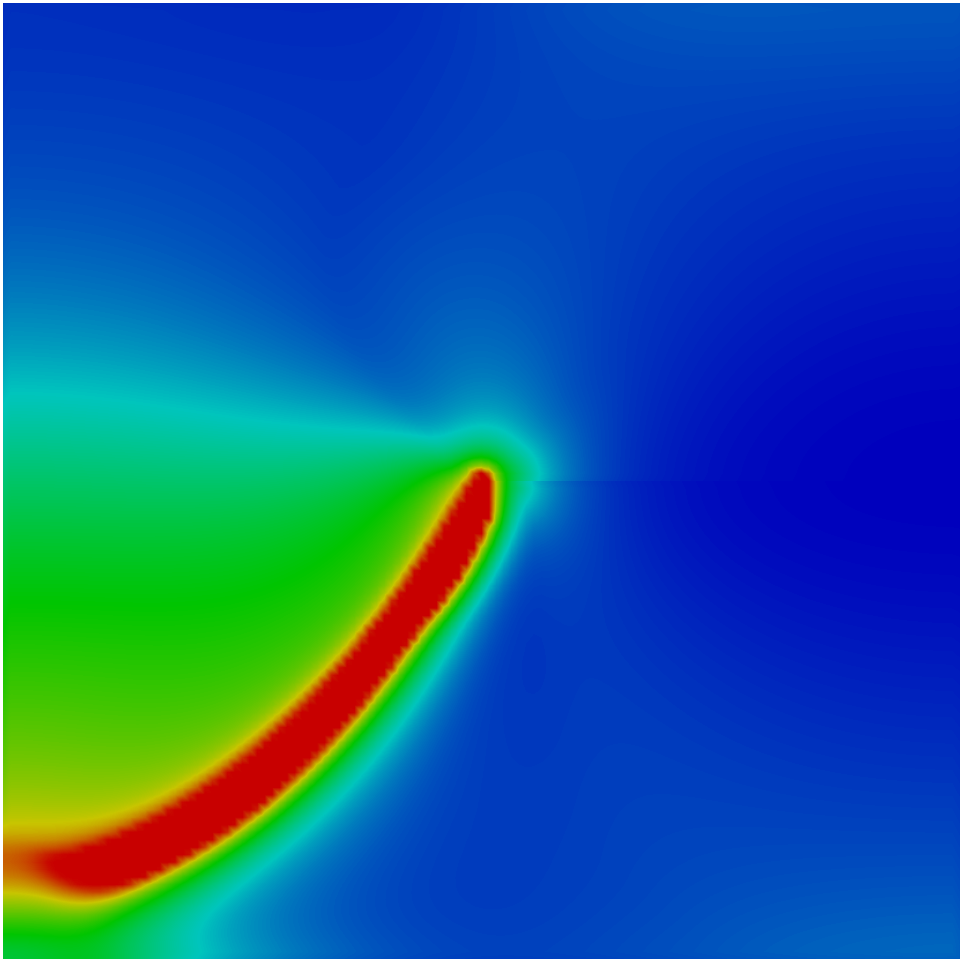}
 \includegraphics[width=0.52\textwidth]{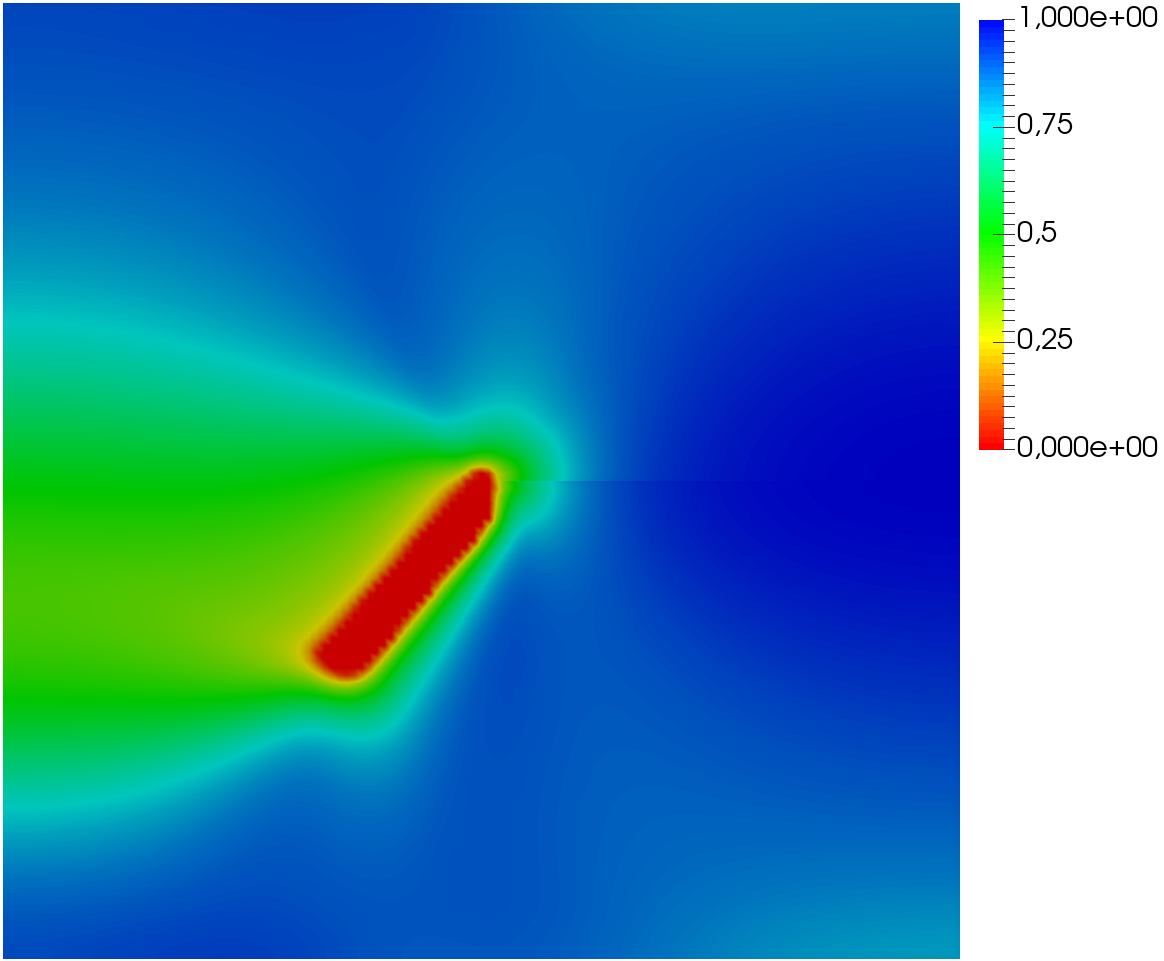}
\end{minipage}
\hfill
\begin{minipage}{0.48\textwidth}
 \includegraphics[width=0.432\textwidth]{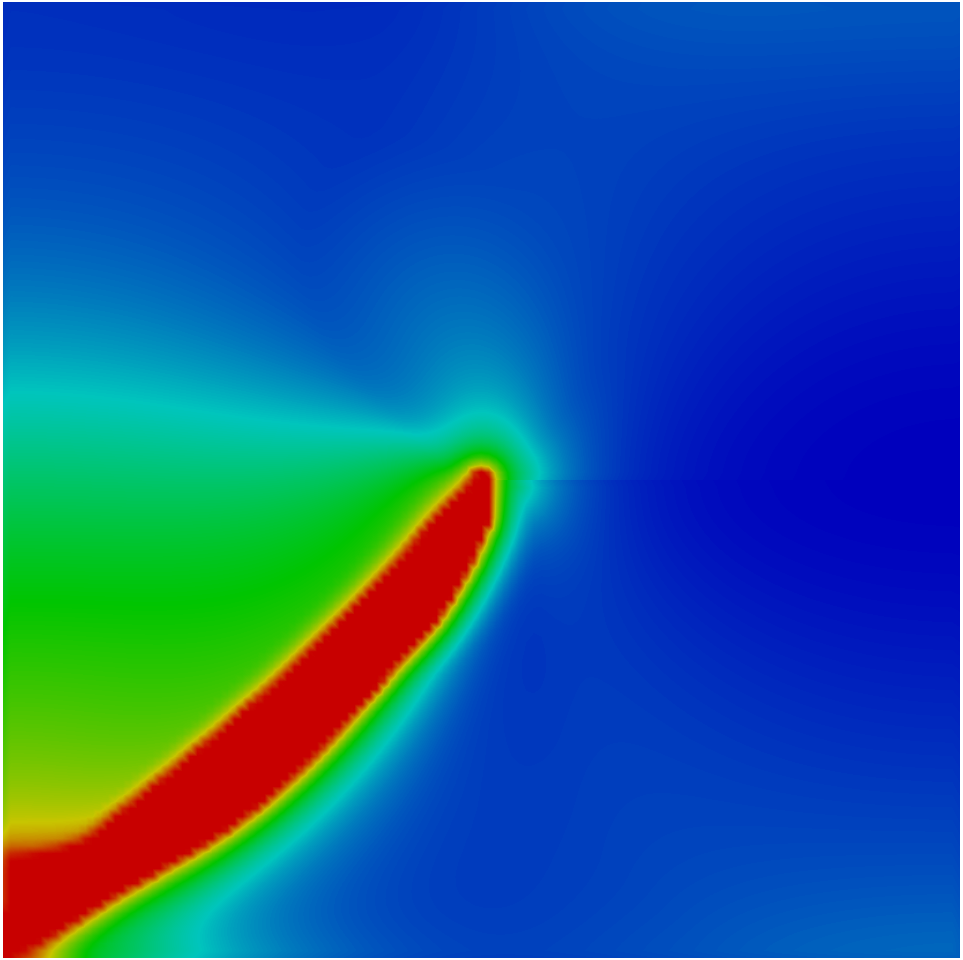}
 \includegraphics[width=0.52\textwidth]{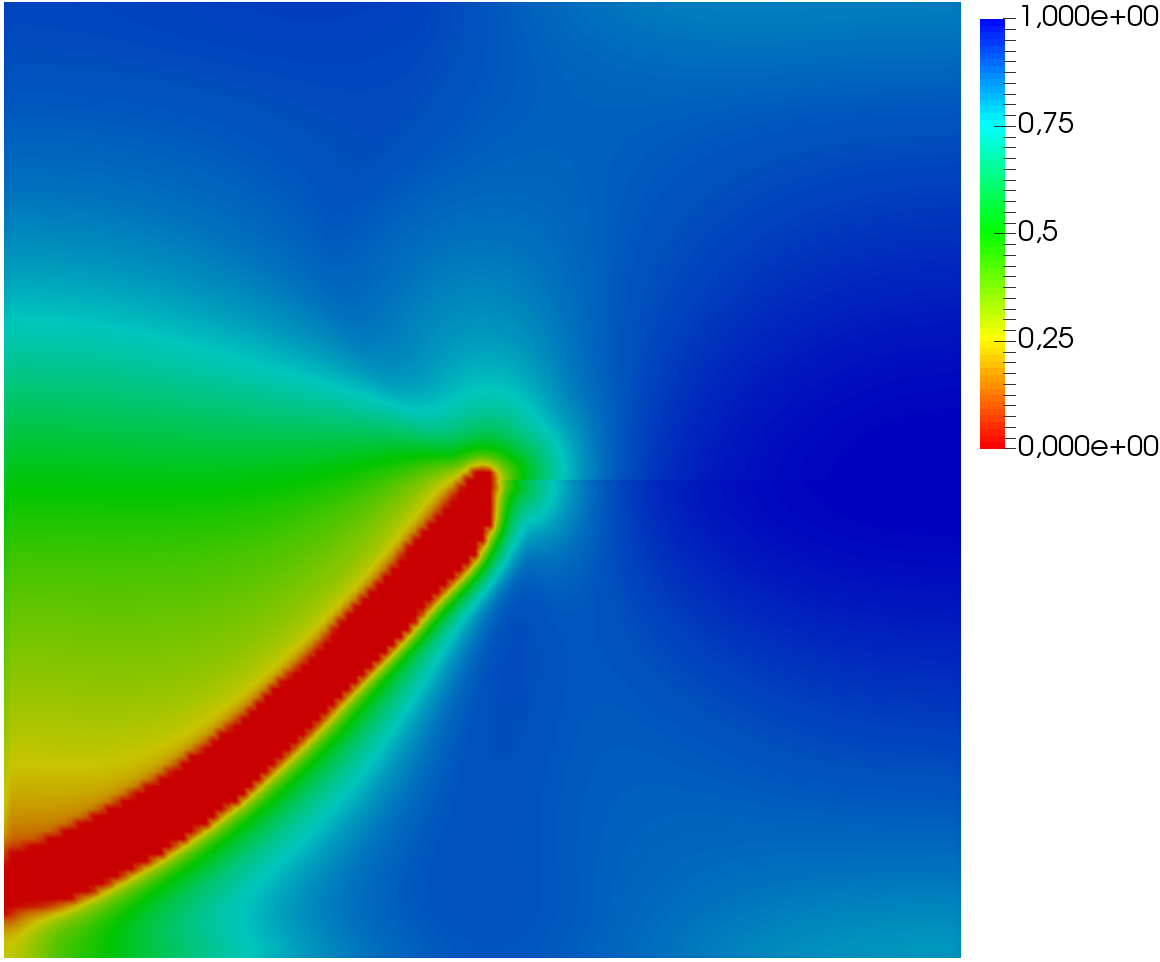}
 \end{minipage}
 \caption{Phase-field function with loading $u_y = 0.033$ and $u_y= 0.042 \si{mm}$ from top to bottom line and for $\nu=0.499$ (left) and for $\nu=0.4999$ (right) with $6$ steps of uniform refinement. }\label{shear_nu_later}
\end{figure}

\begin{figure}[htbp!]
\begin{minipage}{0.48\textwidth}
 \includegraphics[width=0.432\textwidth]{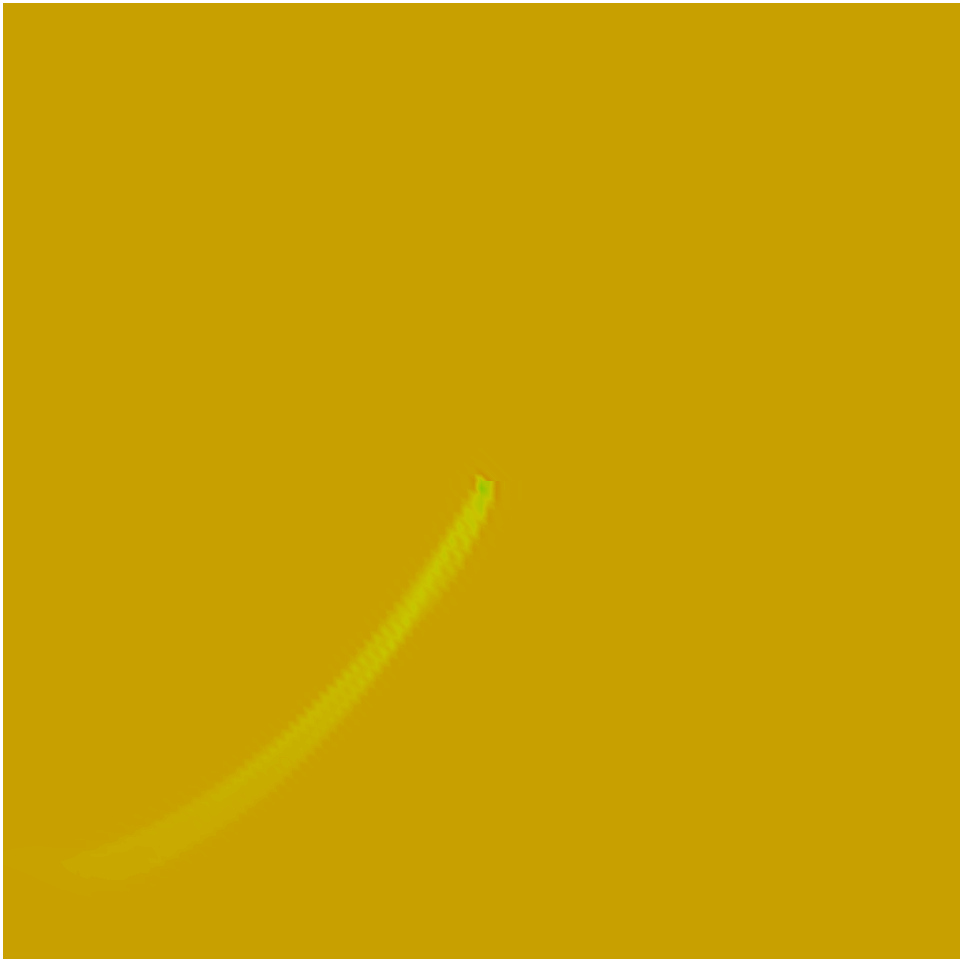}
 \includegraphics[width=0.528\textwidth]{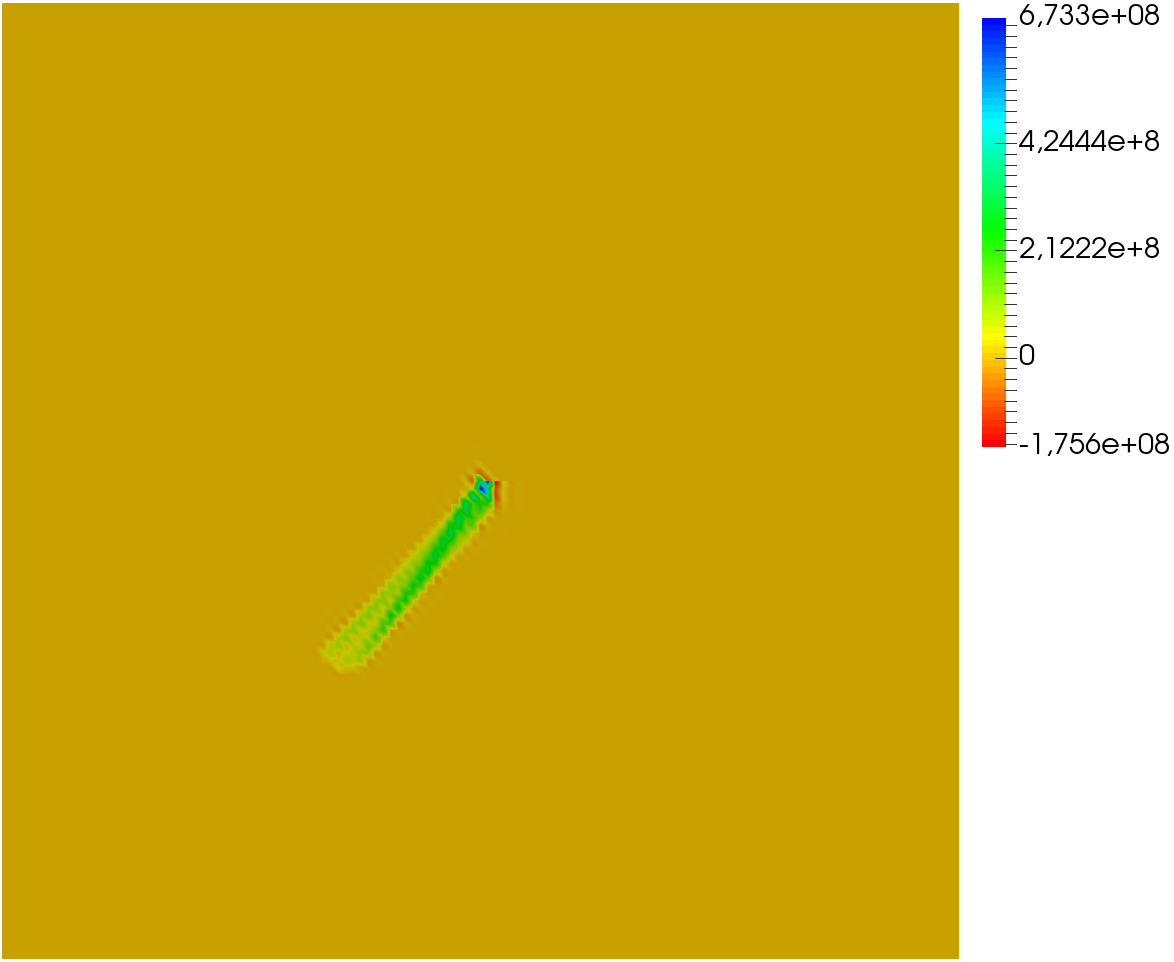}
\end{minipage}
\hfill
\begin{minipage}{0.48\textwidth}
 \includegraphics[width=0.432\textwidth]{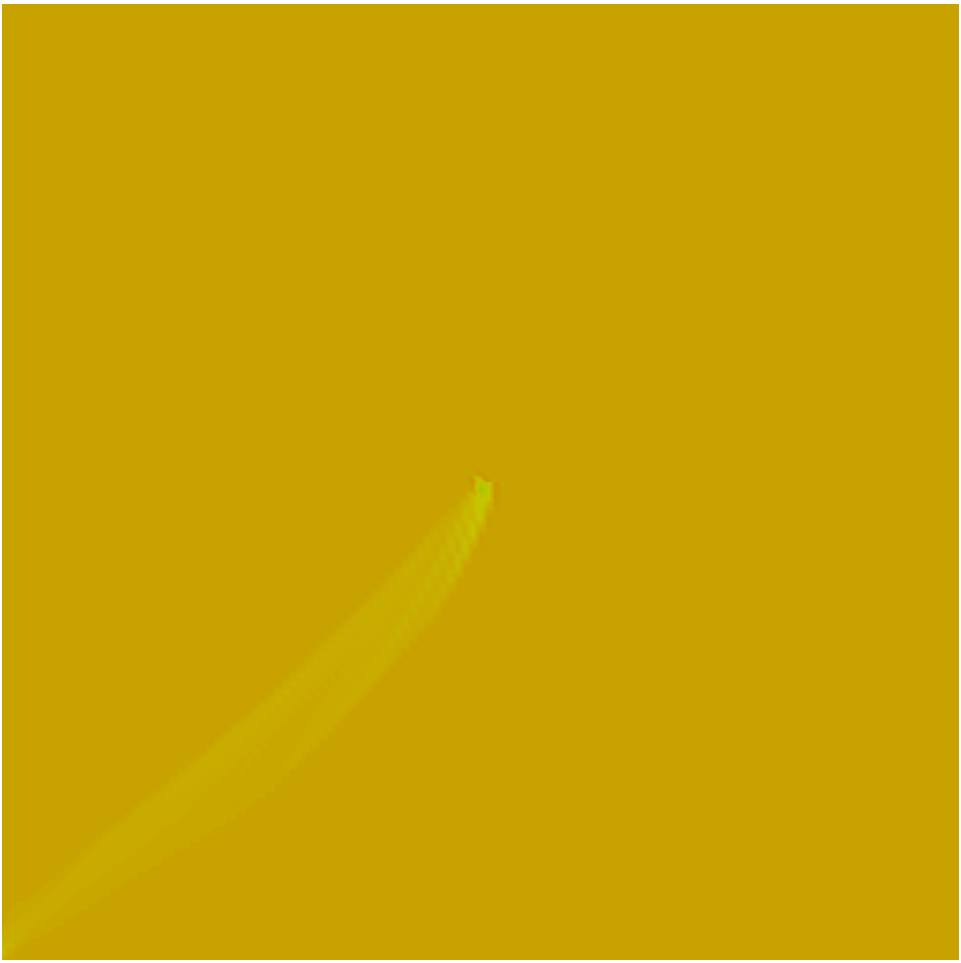}
 \includegraphics[width=0.528\textwidth]{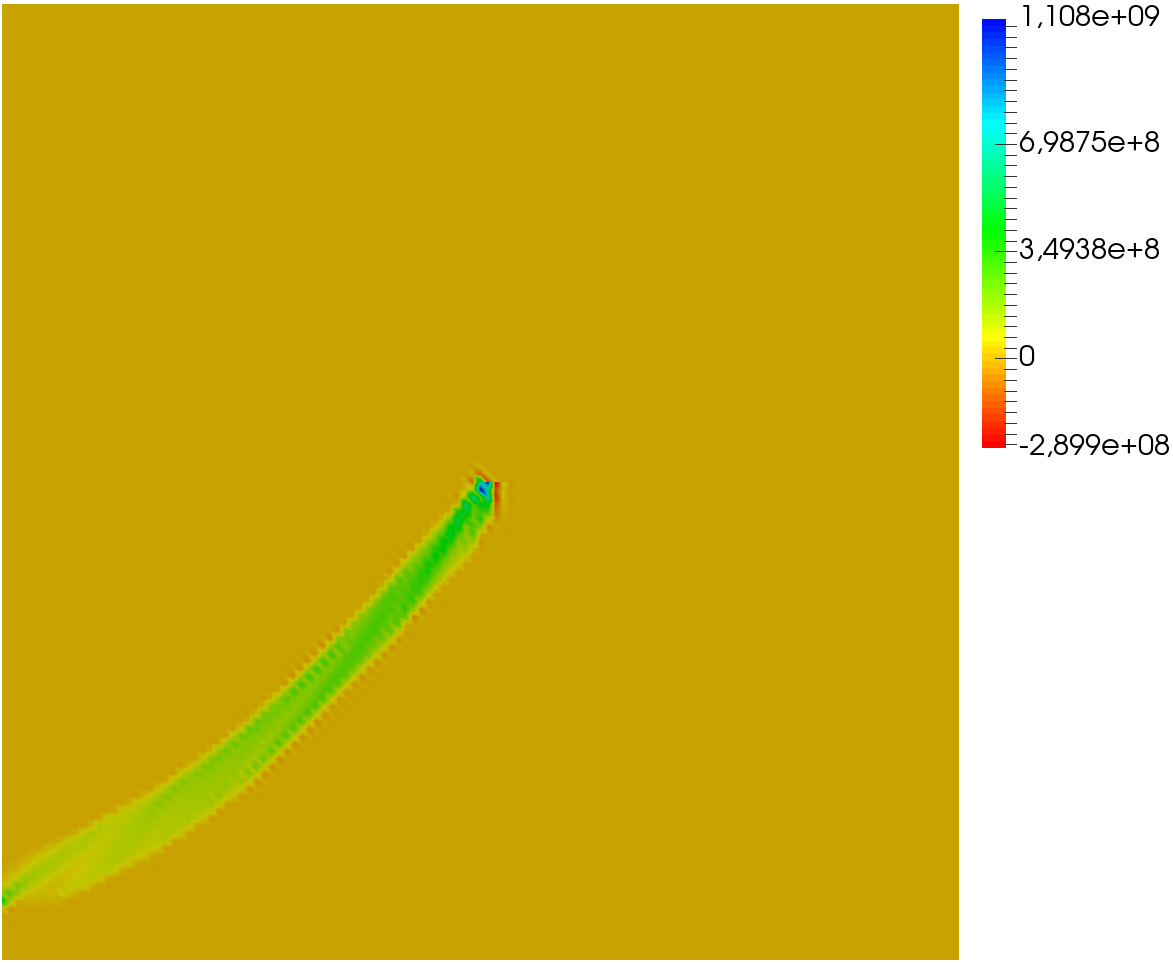}
 \end{minipage}
 \caption{Pressure field with loading $u_y = 0.033$ and $u_y= 0.042 \si{mm}$ from top to bottom line and for $\nu=0.499$ (left) and for $\nu=0.4999$ (right) with $6$ steps of uniform refinement. }\label{shear_press_later}
\end{figure}

Figure \ref{shear_nu_screenshots} presents plots of the phase-field function at certain time steps with an increasing Poisson ratio ($\nu=0.3,\nu=0.49$ and $\nu=0.4999$ from the left to the right). 
The propagation of the crack starts later with an increasing Lam\'{e} coefficient $\lambda$.
For $\nu=0.4999$, the loading of $0.03 \si{mm}$ ($0.03\si{s}$)  does not suffice that the crack tends to the bottom left corner. For this reason, in Figure \ref{shear_nu_later} the phase-field function at later time steps is depicted. Further, plots of the corresponding 
pressure field in Figure \ref{shear_press_later} allow to observe the incompressible behavior with $\nu=0.4999$. The pressure field is plotted at later time steps, where the crack evolves also for high Poisson's ratios. With $\nu=0.4999$ in the right column of Figure \ref{shear_press_later}, 
the necessary pressure for crack propagation is much higher than in less incompressible materials.

\subsubsection{Results of the L-shaped panel test}

In this section, the L-shaped panel setup is tested with higher values for $\nu$, which are listed in Table \ref{l_shaped_nu_table}.
\begin{table}[htbp!]
\centering
\renewcommand*{\arraystretch}{1.4}
\begin{tabular}{|l|l|r|}\hline
\multicolumn{1}{|c}{$\nu$} & \multicolumn{1}{|c}{$\mu$} & \multicolumn{1}{|c|}{$\lambda$} \\ \hline \hline
$0.18$  & $10.95\cdot 10^3$ & $6.18\cdot 10^3$\\ \hline
$0.3$  & $10.95\cdot 10^3$ & $15.88\cdot 10^3$ \\ \hline
$0.4$ & $10.95\cdot 10^3$ & $42.36\cdot 10^3$ \\ \hline
$0.45$ & $10.95\cdot 10^3$ & $95.31\cdot 10^3$ \\ \hline
$0.49$ & $10.95\cdot 10^3$ & $518.91\cdot 10^3$ \\ \hline
$0.499$ & $10.95\cdot 10^3$ & $5464.05\cdot 10^3$ \\ \hline
$0.4999$ & $10.95\cdot 10^3$ & $54739.10\cdot 10^3$ \\ \hline
 \end{tabular}
\caption{Tests with different Poisson ratios approximating $\nu=0.5$ for the L-shaped panel test.}\label{l_shaped_nu_table}
\end{table}

\begin{figure}[htbp!]
\centering
\begin{tikzpicture}[xscale=0.75,yscale=0.75]
\begin{axis}[
    ylabel = Load $F_y$ $\lbrack\si{kN}\rbrack$,
    xlabel = Displacement $\lbrack\si{mm}\rbrack$,
 legend pos=south west, grid =major,
    x post scale = 1.4,
    y post scale = 2.1,
  xtick={-0.3,0,0.3,0.6,0.8,1.0}, 
  ytick={-200,-150,-100,-50,0,25}
  ]
\addplot[green]
table[x=Loading,y=030_4global,col sep=comma] {030_load_displacement_l_shaped.csv}; 
\addlegendentry{$\nu = 0.3$}
\addplot[cyan,dashed]
table[x=Loading,y=040_4global,col sep=comma] {040_load_displacement_l_shaped.csv}; 
\addlegendentry{$\nu = 0.4$}
\addplot[orange,densely dotted]  
table[x=Loading,y=045_4global,col sep=comma] {045_load_displacement_l_shaped.csv};
\addlegendentry{$\nu = 0.45$}
\addplot[blue]  
table[x=Loading,y=049_4global,col sep=comma] {049_load_displacement_l_shaped.csv};
\addlegendentry{$\nu = 0.49$}
\addplot[yellow]  
table[x=Loading,y=0499_4global,col sep=comma] {0499_load_displacement_l_shaped.csv};
\addlegendentry{$\nu = 0.499$}
\addplot[red,dashed]  
table[x=Loading,y=04999_4global,col sep=comma] {04999_load_displacement_l_shaped.csv};
\addlegendentry{$\nu = 0.4999$}
\end{axis}
\end{tikzpicture}
\caption{Load-displacement curves of the L-shaped panel test for different Poisson ratios and $4$ steps of uniform refinement.}\label{l_shaped_nu_plot}
\end{figure}
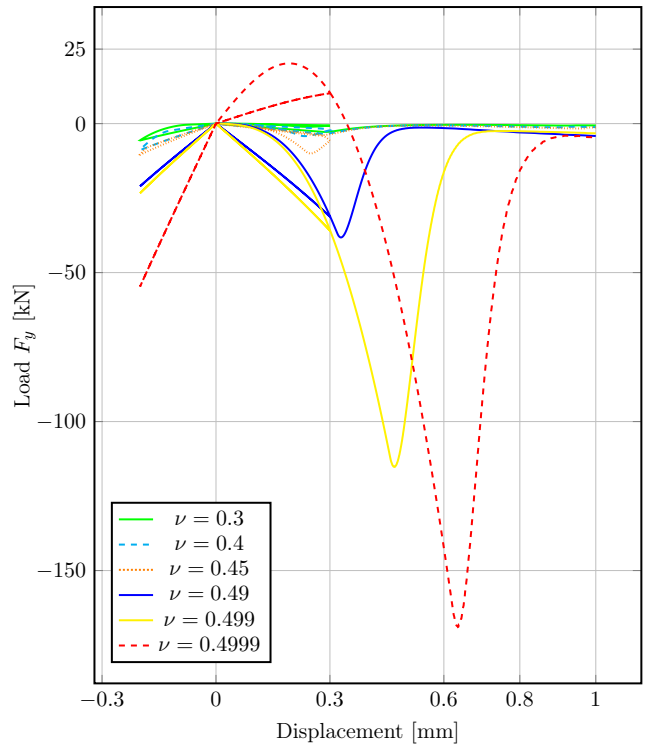

As depicted in Figure \ref{l_shaped_nu_plot}, the crack does not just evolve later in time and with a larger loading force but the whole course of the \selectlanguage{ngerman}load"=displacement\selectlanguage{english} curves changes. This can be observed especially for $\nu=0.4999$ (the red dotted curve).

\begin{figure}[htbp!]
\begin{minipage}{0.48\textwidth}
 \includegraphics[width=0.305\textwidth]{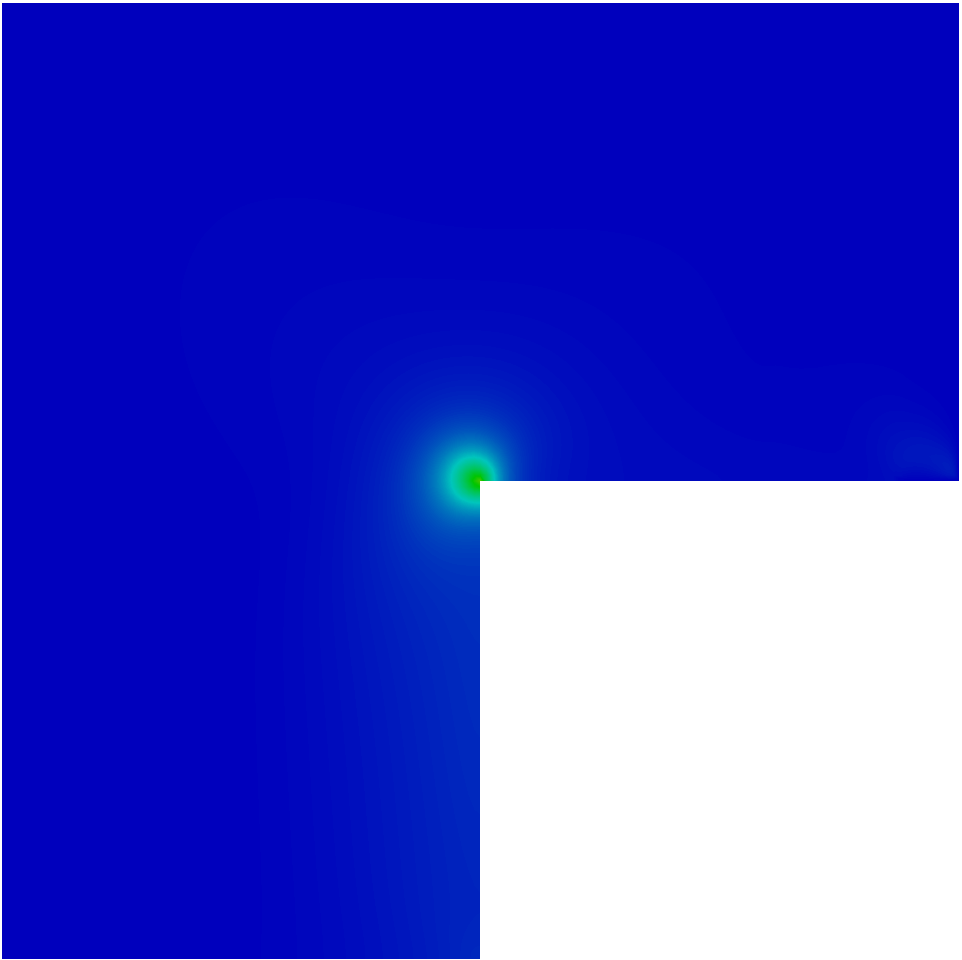}
 \includegraphics[width=0.305\textwidth]{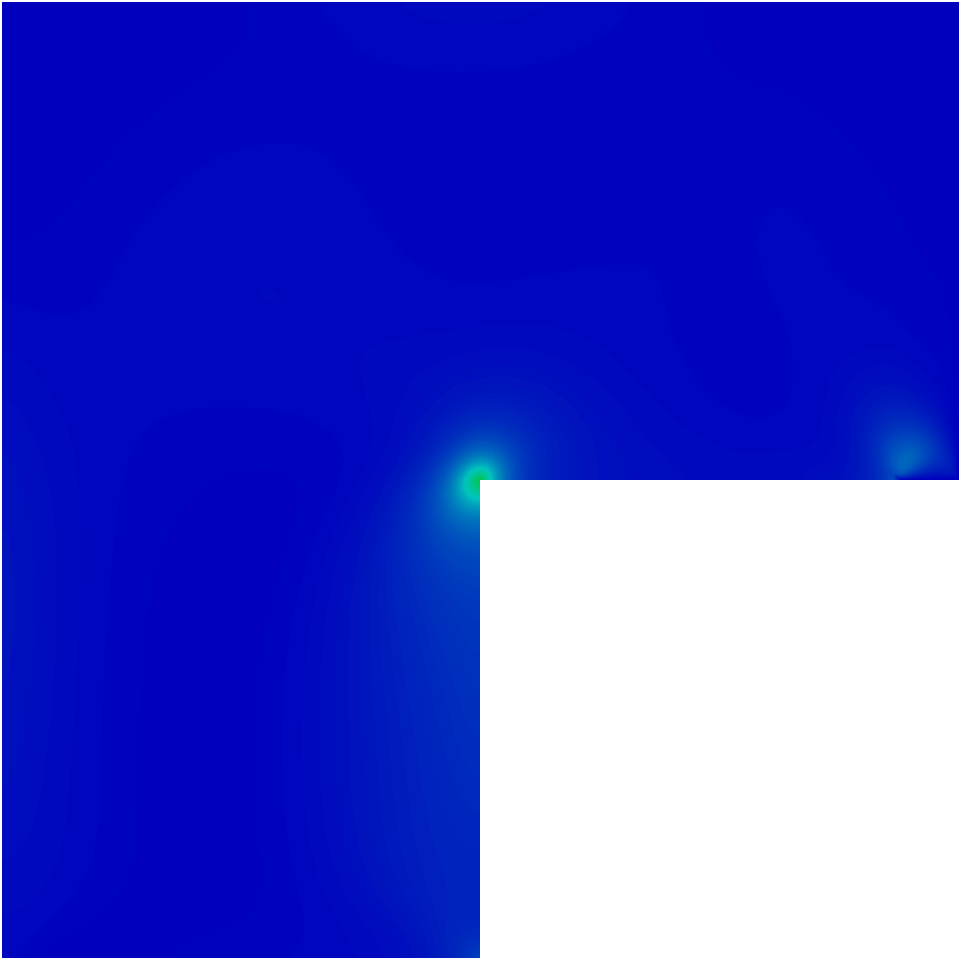}
 \includegraphics[width=0.369\textwidth]{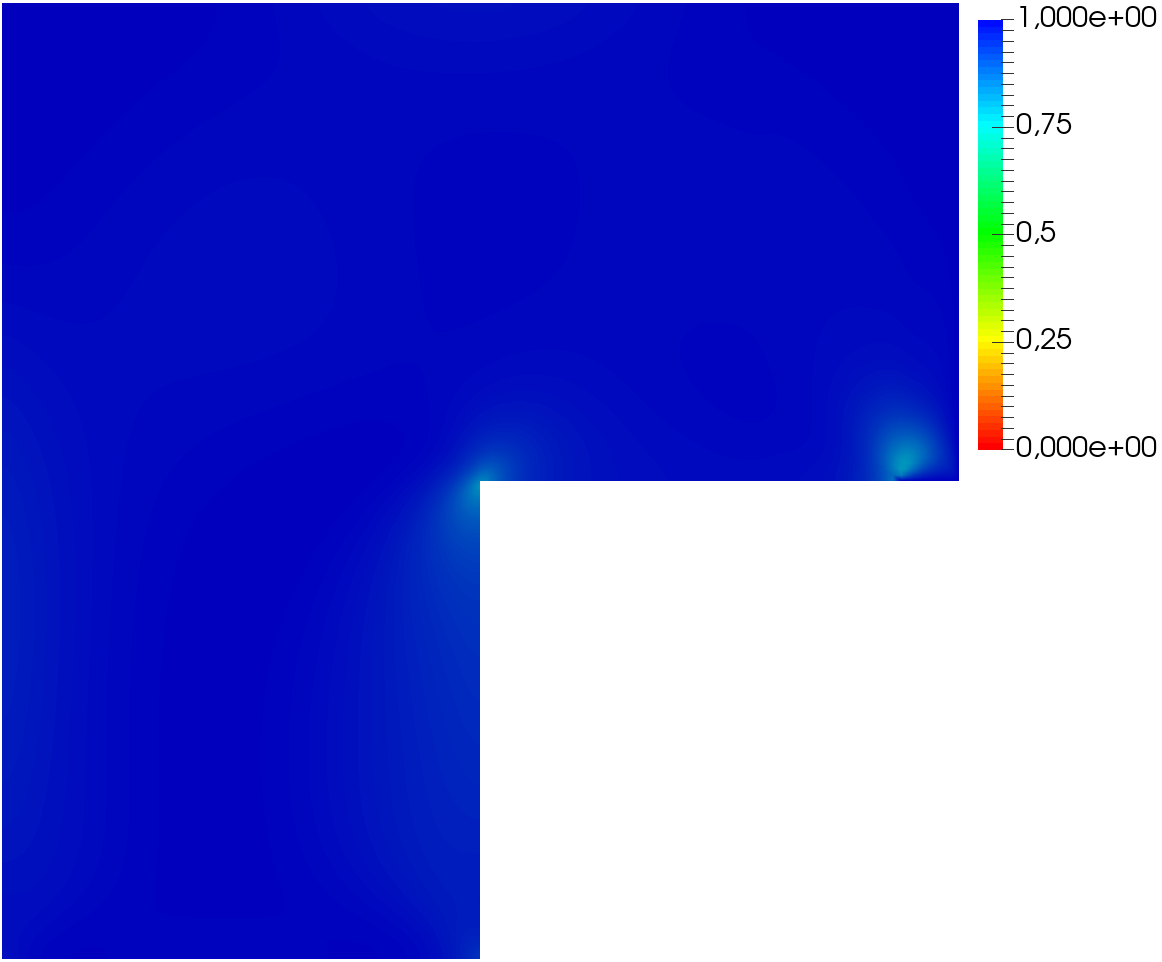}
\end{minipage}
\hfill
\begin{minipage}{0.48\textwidth}
 \includegraphics[width=0.305\textwidth]{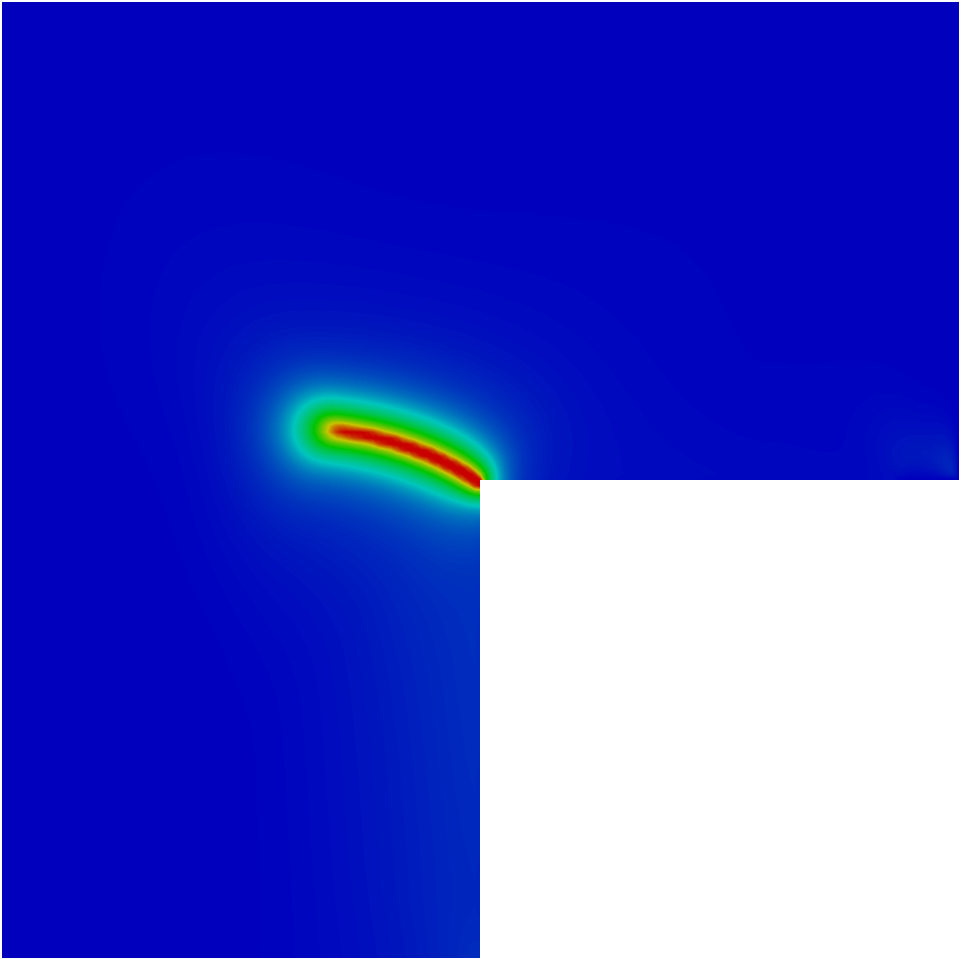}
 \includegraphics[width=0.305\textwidth]{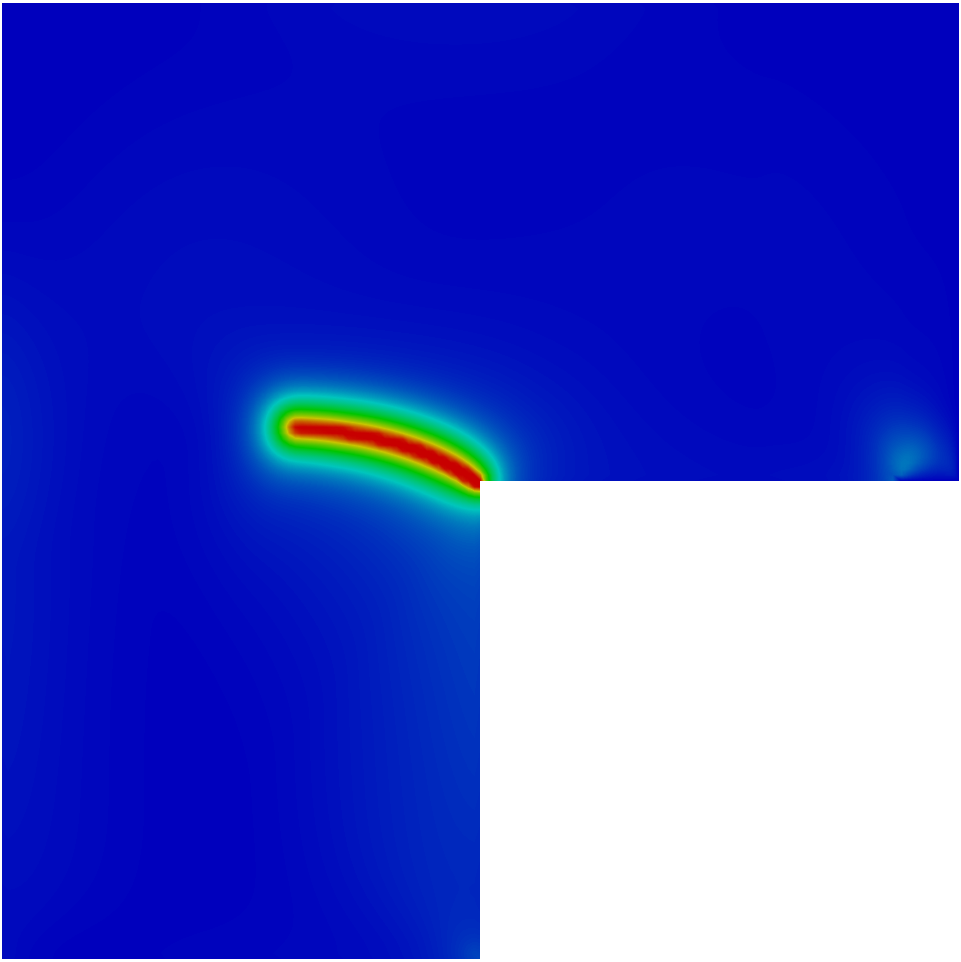}
 \includegraphics[width=0.369\textwidth]{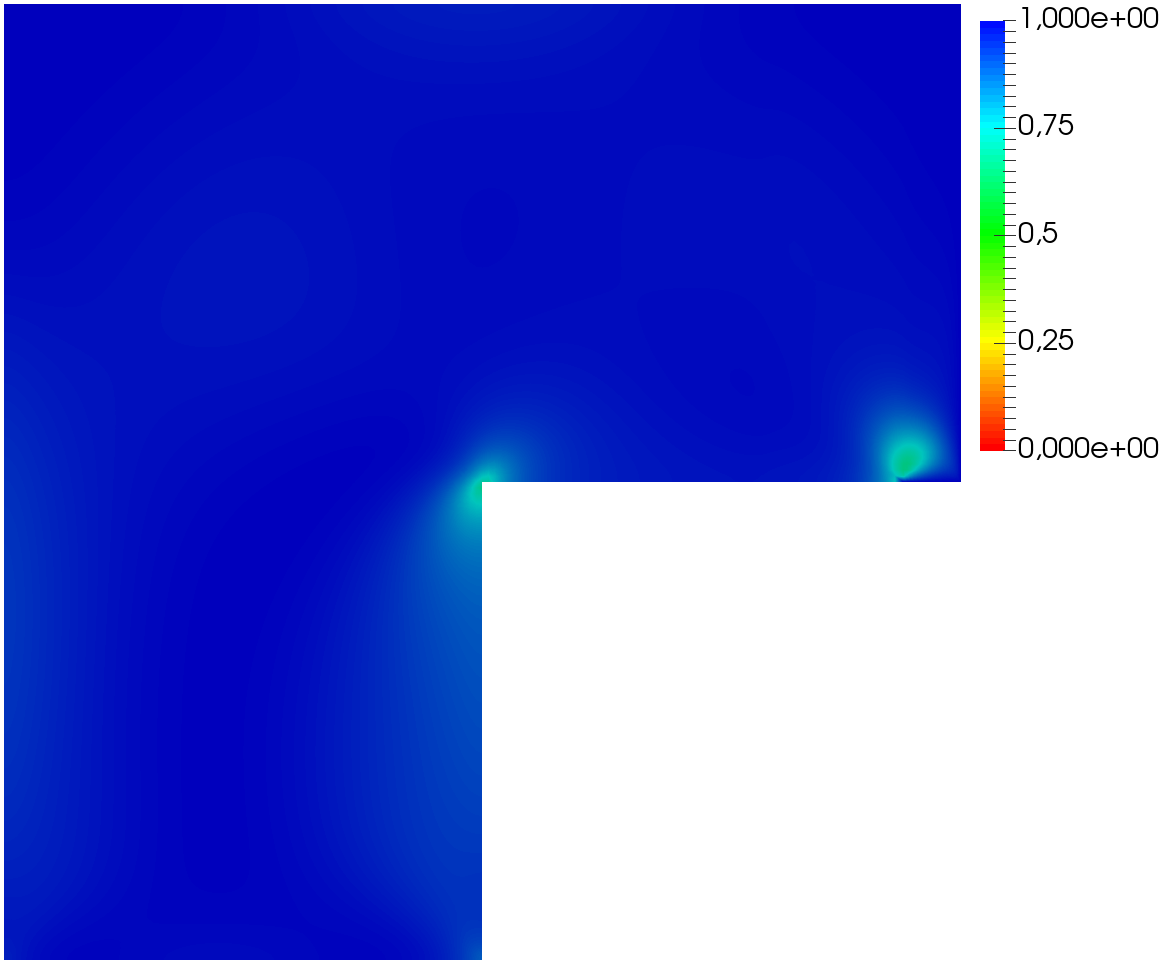}
 \end{minipage}
 \hfill
 \begin{minipage}{0.48\textwidth}
 \includegraphics[width=0.305\textwidth]{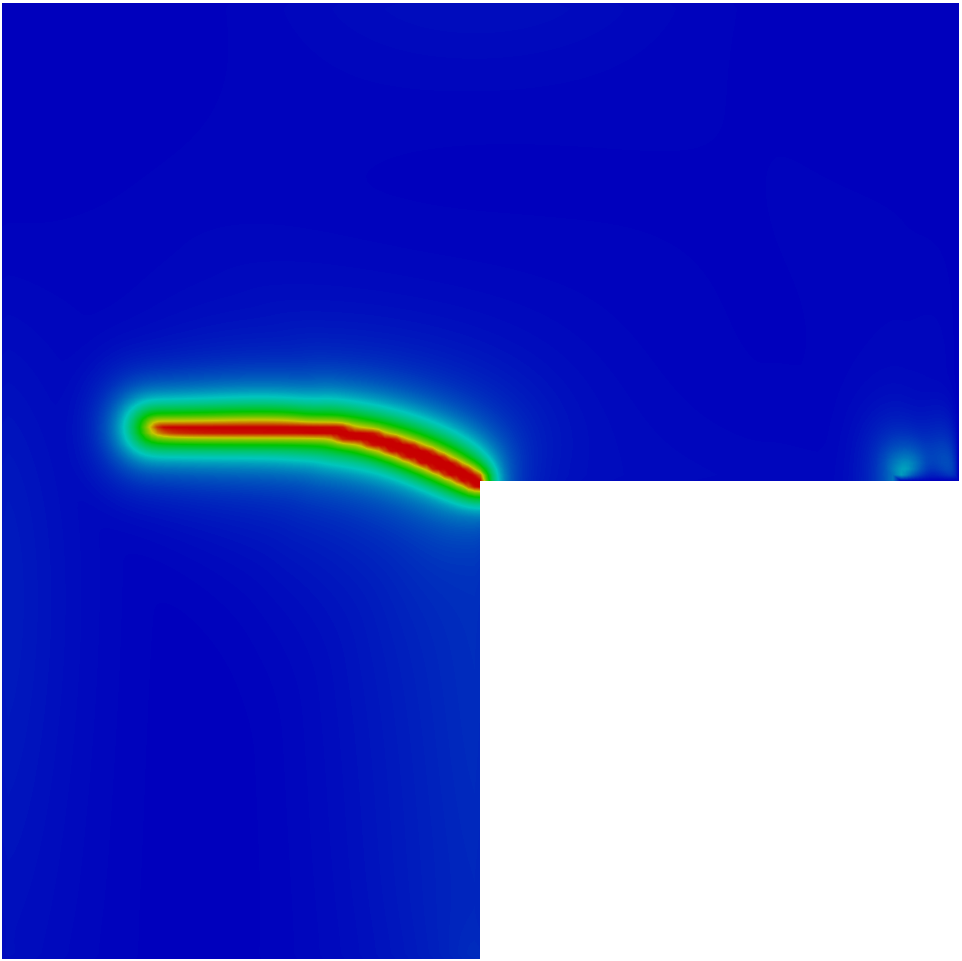}
 \includegraphics[width=0.305\textwidth]{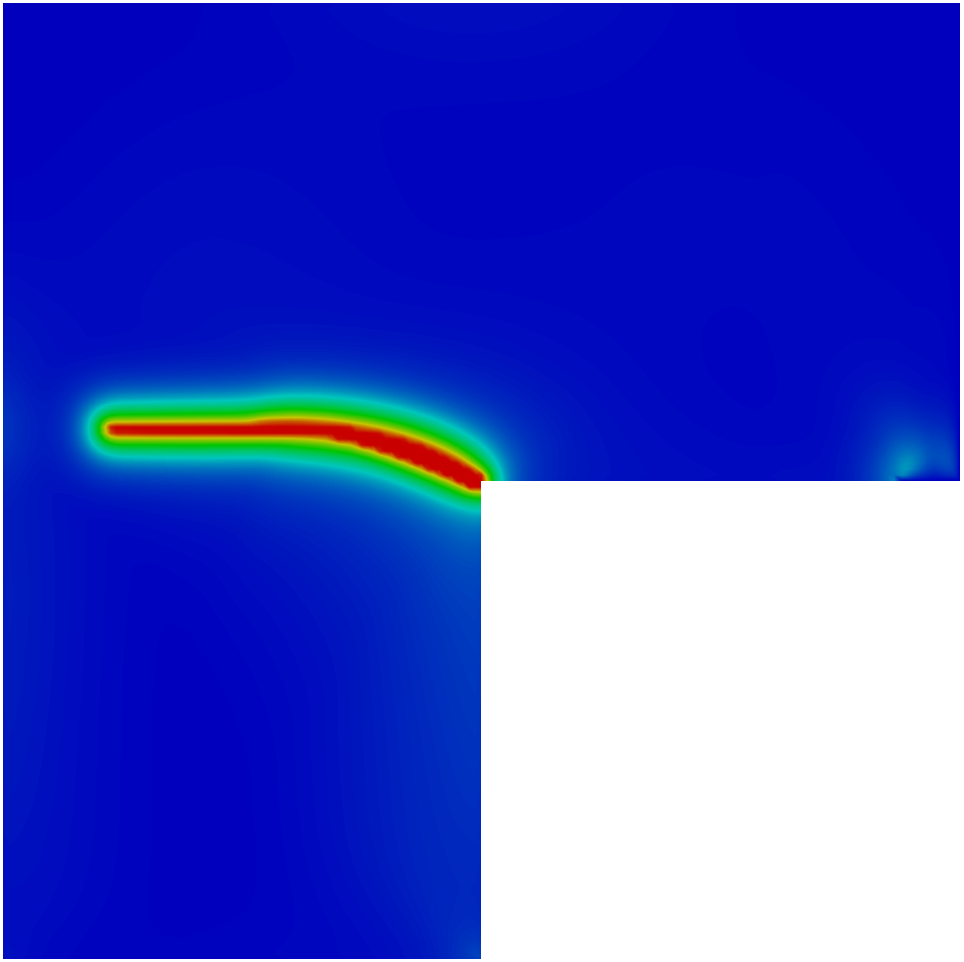}
 \includegraphics[width=0.369\textwidth]{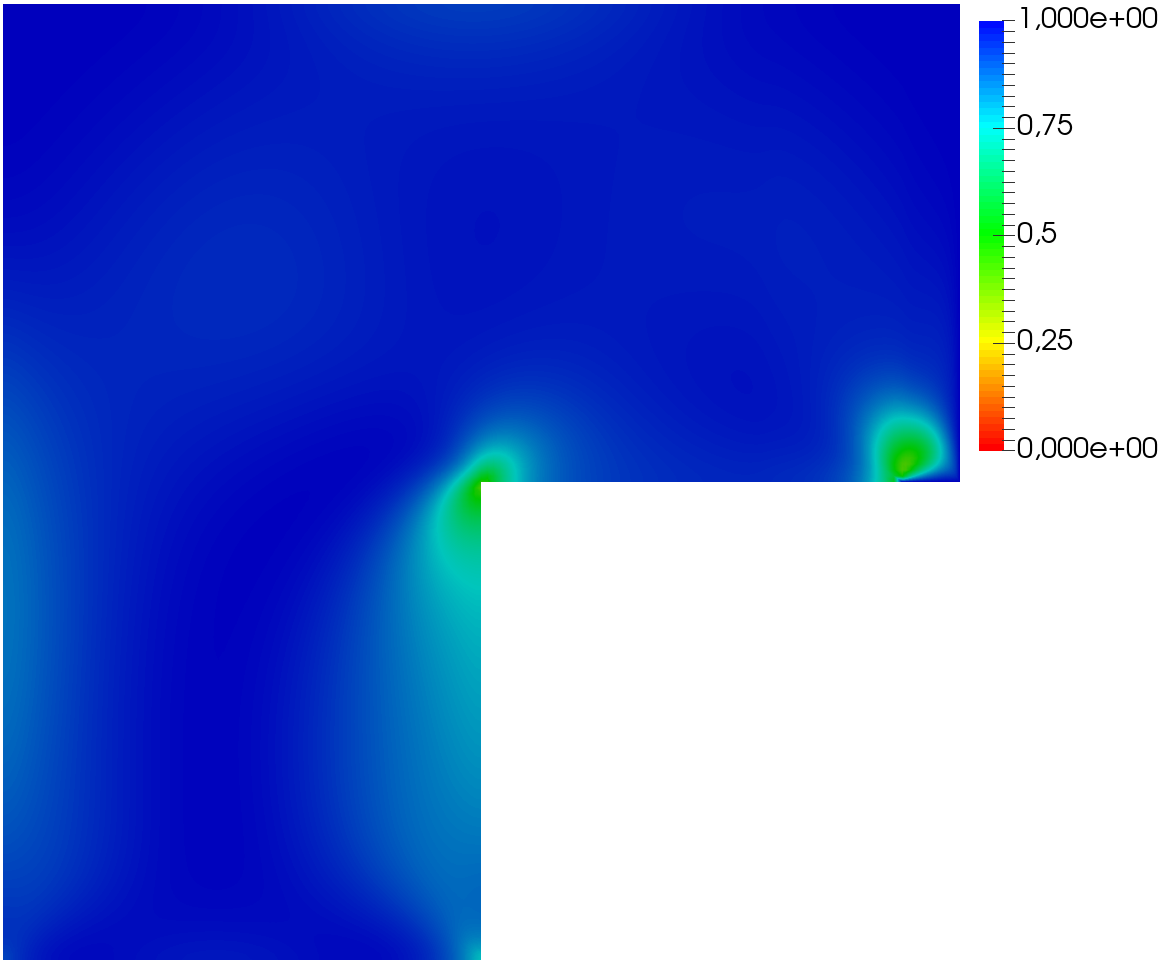}
\end{minipage}
\hfill
\begin{minipage}{0.48\textwidth}
 \includegraphics[width=0.305\textwidth]{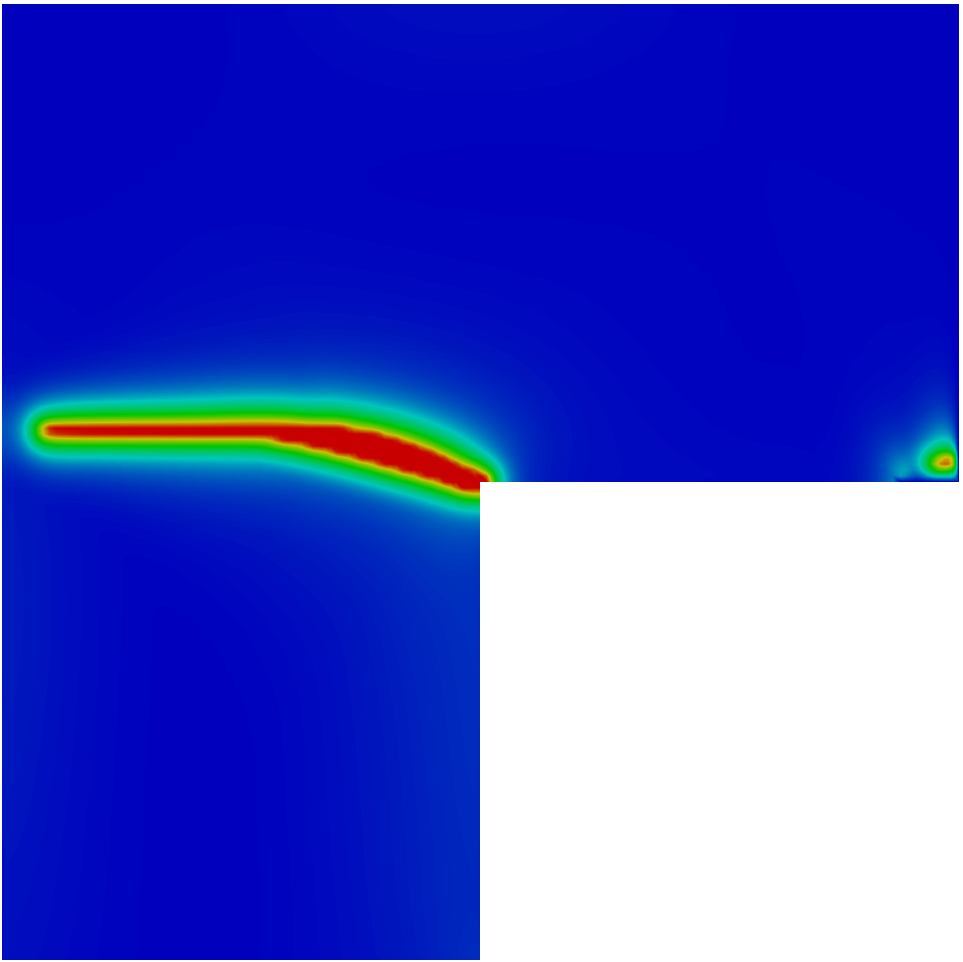}
 \includegraphics[width=0.305\textwidth]{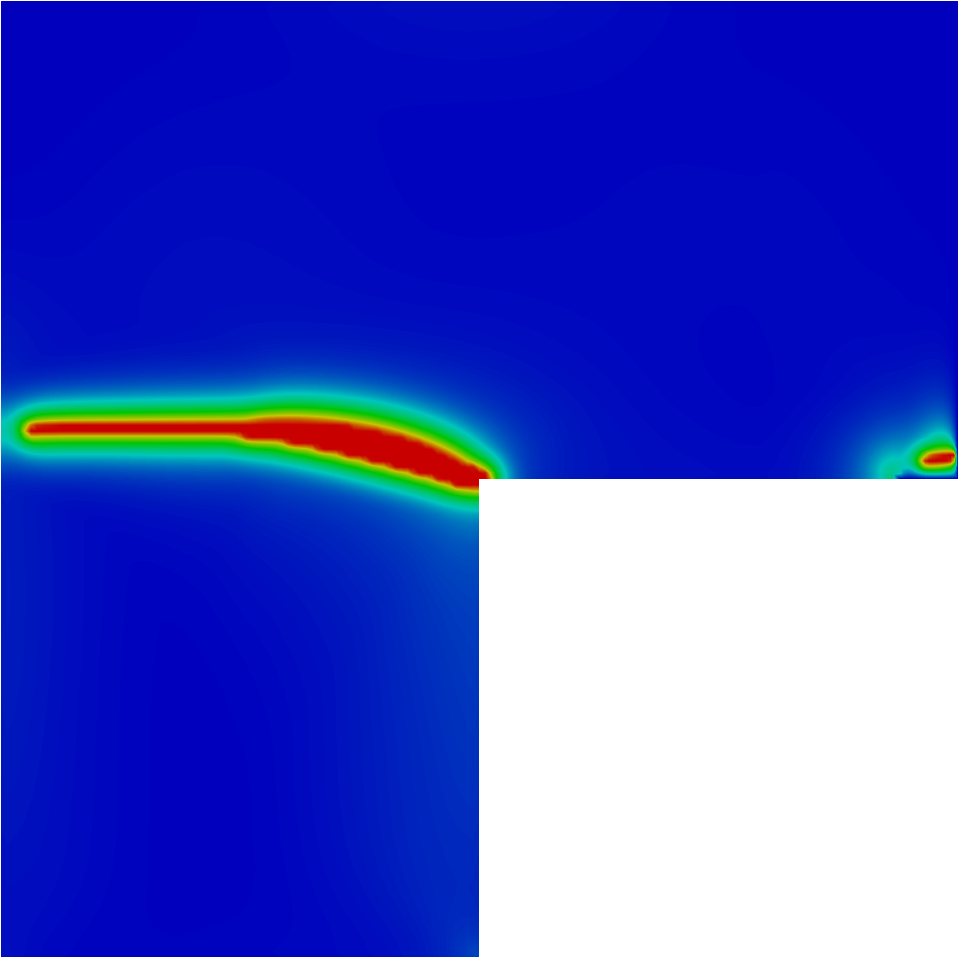}
 \includegraphics[width=0.369\textwidth]{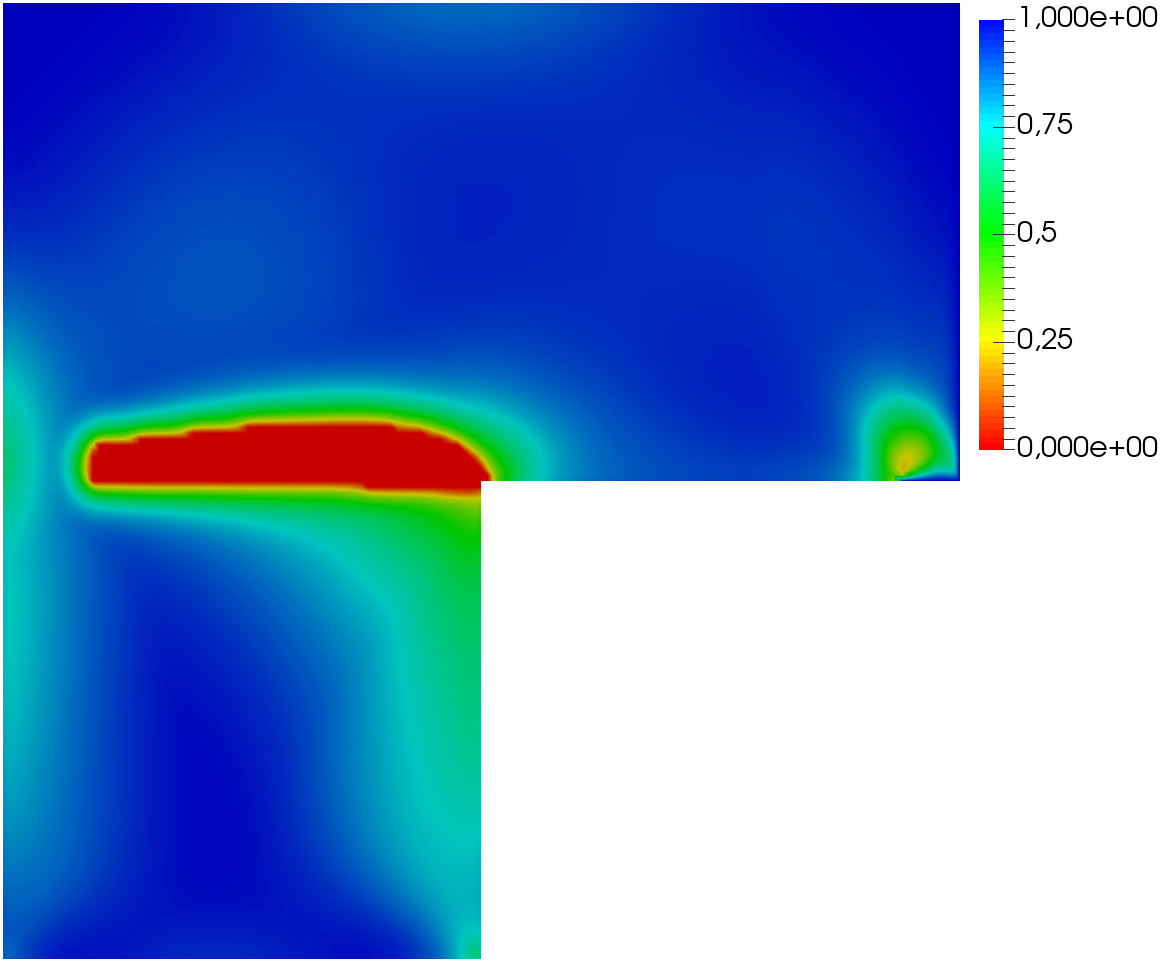}
 \end{minipage}
 \caption{Phase-field function with loading $u_y = 0.22, 0.3, 0.45$ and $1.0 \si{mm}$ from top to bottom line and for $\nu=0.18$ (left), for $\nu=0.45$ (middle) and for $\nu=0.4999$ (right) with $4$ steps of refinement. }\label{l_shaped_nu_screenshots}
\end{figure}

Plots of the phase-field function at certain time steps with increasing Poisson ratios ($\nu=0.3,\nu=0.45$ and $\nu=0.4999$ from the left to the right) are depicted in Figure \ref{l_shaped_nu_screenshots}. 
The propagation of the crack starts later in time with an increasing Lam\'{e} coefficient $\lambda$. It stands out that the crack seems to grow wider not just in the corner at the midpoint with an increasing Poisson's ratio $\nu$.


\section{Conclusions}\label{Conclusion}
The focus of this work was to 
develop a phase-field model for fractures in incompressible materials.
To ensure stability, we derived a mixed system of a standard
phase-field model. 
As it is well-known for mixed systems with inf-sup stability, the corresponding 
finite element spaces have to be chosen carefully. 
We use biquadratic elements for the displacement function and bilinear shape functions for the hydro-static pressure variable.
Detailed discussions of the mixed formulation for this phase-field fracture problem
were provided in Section 3 and 4.
In Section 5 we adopted the settings of 
two well-known numerical test and designed 
a series of numerical studies. The main goals were 
a comparison of different finite element orders for the standard phase-field
model in order to study the influence of higher-order finite elements 
on phase-field modeling. Here, we observed 
small changes for the single edged notched shear test,
but significant changes for the L-shaped panel test. 
This allows assuming, 
that the L-shaped panel test is more sensitive with respect to the choice of finite elements. 
Then, we conducted studies on meshes with different levels of uniform refinement and proposed tests with different Poisson ratios $\nu$ approximating the incompressible limit $\nu =0.5$. 
The \selectlanguage{ngerman}load"=displacement\selectlanguage{english} curves of both tests show a correlation between an increasing Poisson ratio and
a stronger loading force. 
Specifically, for increasing Poisson's ratios higher stresses 
are observed before cracking. 
Future work is to extend this model to an a posteriori error estimation 
and adaptive refinement strategies.



\section*{Acknowledgments}
This work has been supported by the German Research Foundation, Priority Program 1748 (DFG SPP 1748) named
\textit{Reliable Simulation Techniques in Solid Mechanics. Development of
Non-standard Discretization Methods, Mechanical and Mathematical Analysis}. Our subproject within the SPP1748 reads \textit{Structure Preserving Adaptive Enriched Galerkin Methods for Pressure-Driven 3D Fracture Phase-Field Models} (WI 4367/2-1 and WO 1936/5-1).


\end{document}